\documentclass[]{elsarticle}
\usepackage{amsmath} 
\usepackage{mathtools}
\usepackage{eurosym}
\usepackage{amssymb}
\usepackage[printonlyused]{acronym}
\usepackage[labelformat = simple,format=plain,bf,font={footnotesize},margin={0.05\linewidth,0.05\linewidth},skip=8pt,tableposition=top,figureposition=bottom,singlelinecheck=false]{caption}
\usepackage[table,dvipsnames]{xcolor}
\usepackage{booktabs}
\usepackage{calc}
\usepackage[]{units}
\usepackage{array}
\usepackage{longtable}
\usepackage{tabu}
\usepackage{slashbox}
\usepackage{setspace}
\usepackage[super]{nth}
\usepackage[hidelinks]{hyperref}
\usepackage[hyphenbreaks]{breakurl}
\usepackage{P3defs}

\journal{European Journal of Operational Research}

\begin{document}

\begin{frontmatter}

\title{Calibration algorithm for spatial partial equilibrium models with conjectural variations}

\author[IED,IfA]{Tobias Baltensperger\corref{correspondingauthor}}
\cortext[correspondingauthor]{Corresponding author. Tel.: +41 44 632 44 73;}
\ead{t.baltensperger@usys.ethz.ch}

\author[ZHAW]{Rudolf M. F\"{u}chslin}
\author[IED]{Pius Kr\"{u}tli}
\author[IfA]{John Lygeros}

\address[IED]{Institute for Environmental Decisions (IED), ETH Z\"{u}rich, 8092 Z\"{u}rich, Switzerland}
\address[IfA]{Automatic Control Laboratory (IFA), ETH Z\"{u}rich, 8092 Z\"{u}rich, Switzerland}
\address[ZHAW]{Institute of Applied Mathematics and Physics (IAMP), ZHAW Z\"{u}rich University of Applied Sciences, 8401 Winterthur, Switzerland}

\begin{abstract}
When calibrating spatial partial equilibrium models with conjectural variations, some modelers fit the suppliers' sales to the available data in addition to total consumption and price levels. While this certainly enhances the quality of the calibration, it makes it difficult to accommodate user-imposed bounds on the model parameters such as restricting the market power parameters to the interval [0,1], which is a common requirement in conjectural variations approaches. We propose an algorithm to calibrate the suppliers' sales and simultaneously deal with user-defined bounds on parameters.
To this end, we fix the suppliers' sales at reference values and obtain the marginal costs for each supplier and market. We then limit the market power parameters to the interval [0,1], and calculate intervals of anchor prices and price elasticities that reproduce the reference supplier sales in the state of equilibrium. If these intervals do not contain the reference price elasticities and prices, we face a mismatch between reality and the model mechanics. We resolve this by altering the reference sales for the critical suppliers, and iterate. Thereby, the user controls whether price elasticities and anchor prices should be close to their reference values, or the suppliers' sales.
The algorithm is tested on data from the European gas market, and required less than one minute to identify calibrated parameters. Our algorithm is widely applicable, since it is based on mild (and common) underlying assumptions and can be configured to suit a specific purpose thanks to the inclusion of user-defined bounds on all relevant parameters.

\end{abstract}


\end{frontmatter}

\section{Introduction} \label{sec:P3_Introduction}

With the advancement in computer technology, an increasing number of market studies is based on mathematical models. A popular way of modeling a network of markets is via a spatial partial equilibrium model including \ac{CV}. However, the calibration problem associated with these models turns out to be a \ac{MPEC}, as one aims to minimize the deviations between some of the parameters of the model and their reference values subject to the original equilibrium problem. Unfortunately, \ac{MPEC}s have a non-convex solution space and are thus difficult and time-consuming to solve for large models. 

An alternative method is used by \citet{Lise2008}, \citet{Cobanli2014} and other gas market modelers. First, consumption per node and time period is fixed and, under the assumption of perfectly competitive markets, the nodal prices are computed. Then the inverse demand functions are derived from predefined price elasticities and the obtained consumption-price pairs, and the equilibrium is computed again under the assumption of imperfectly competitive markets; this reduces consumption and increases price levels. The inverse demand functions are successively shifted outwards until the original consumption levels are matched. 

The resulting nodal prices also depend on the predefined market power parameters (conjectural parameters) of the suppliers and the price elasticities of the consumers; some modelers consider these parameters as well to achieve better results. \citet{Chyong2014e} tune the price elasticities while leaving the market power parameters at their fixed values. While the achieved equilibrium is impressively close to what is observed in reality, this method can yield unrealisticly price-inelastic consumers.

In contrast, \citet{Garcia-Alcalde2002}, \citet{Liu2006} and other electricity market modelers \citep{LopezdeHaro2007,Liu2007,Diaz2010,Kaminski2011,Diaz2011,Diaz2014} tune the market power parameters while leaving price elasticities fixed. The strength of this approach is that the information contained in the level and spatial distribution of the sales of the suppliers is used for calibration, instead of basing the values of the market power parameters on the experience of the modeler. 
On the downside, the marginal supply costs of each trader have to be known. While this can be approximated by the production costs in the single-market settings of \citet{Garcia-Alcalde2002} and the other modelers, the task is more difficult in a network of markets, because transport costs, congestion fees, etc.\ have to be included. Moreover, the market power parameters are allowed to take arbitrary values, even though, in the absence of cartels, only values between 0 and 1 can be associated with the \ac{CV} approach, see for instance \citet[Chapter 12]{Tremblay2012}. 

In this work, we propose a new calibration algorithm that bridges these gaps by finding market power parameters in the interval $[0,1]$ based on reference data on sales from individual suppliers to consumers, while maintaining wholesale price levels, price elasticities, and nodal consumption within user-defined bounds.

\section{Calibration Framework} \label{sec:P3_CalibrationFramework}

\subsection{Oligopolistic market representation}
We describe an oligopolistic market by an equilibrium problem comprising the \ac{KKT} conditions of the optimization problems of the $\bar{f}$ traders supplying the market, and the market clearing condition.
\begin{align}
0 &\leq\, -\lambda -\MP_f \frac{d\IDF(s)}{ds} q_{f} + \phi_{f} \perp q_{f} \,\geq 0 \quad \forall f \in \F=\{f_1, \ldots, f_{\bar{f}}\}. \label{eqn:P3_dqC} \\
0 &\leq\,- \lambda + \IDF(s) \perp \lambda \geq 0 \label{eqn:P3_dpiC}
\end{align}

$\lambda$ is the market price, $\MP_{f}$ is the market power parameter of trader $f$ in the set of all traders $\F$ with access to market, $\IDF(s)$ is the inverse demand function, $s:=\sum \limits_{f \in \F} q_{f}$ is the total consumption in the market, $q_{f}$ is the quantity sold by trader $f$, and $\phi_f$ is trader $f$'s marginal cost of supplying $q_f$.

We require the inverse demand function $\IDF(s)$ to be bijective and have the following characteristics in its anchor point $(\sO,\lambdaO)$: 
\begin{align}
\IDF(\sO) &= \lambdaO \geq 0, \label{eqn:P3_lambda0} &
\left. \frac{d\IDF(s)}{d s}\right|_{\sO} &= \frac{\lambdaO}{\sO \cdot \eta} < 0, 
\end{align}
where $\sO > 0$ is the anchor consumption, $\lambdaO \geq 0$ the anchor price, and $\eta:=\left.\frac{\partial s}{\partial \lambda} \right|_{\lambdaO}\cdot \frac{\lambdaO}{\sO} <0$ the price elasticity of demand in $(\sO,\lambdaO)$. While $\sO$, $\lambdaO$, and $\eta$ are subject to calibration, as discussed in the next section, the functional form of $\IDF(s)$ can be freely chosen by the modeler.

\subsection{Problem statement} \label{sec:P3_ProblemStatement}
When calibrating \ac{CV} models, the goal is to tune parameters $\sO$, $\lambdaO$, $\eta$, and $\MP$, such that the resulting model equilibrium ($\lambdaStar$, $\qStar$, $\sStar$, $\phiStar$) matches reference values ($\lambdaRef$, $\qRef$, $\sRef := \sum \limits_{f \in \F} \qRef_{f}$, $\phiRef$). $\MP$ denotes the vector $[\MP_1, \ldots, \MP_{\bar{f}}]^T$ or point $(\MP_1, \ldots, \MP_{\bar{f}})$, depending on context; $q_f$ and $\phi_f$ are abbreviated similarly.

In fact, by defining the market power parameters
\begin{equation} 
\MP_f(\lambdaO,\qRef_f,\sO,\phiRef_f,\eta) := \frac{\lambdaO -\phiRef_f}{\frac{\lambdaO}{\sO \cdot (-\eta)} \cdot \qRef_f} \quad \forall f \in \F, \label{eqn:P3_MPCalibChoice}
\end{equation}
and setting the anchor parameters $\sO=\sRef$, $\lambdaO=\lambdaRef$, and $\eta = \etaRef$, we can force the equilibrium to be in ($\lambdaRef$, $\qRef$, $\sRef$, $\phiRef$) and thereby successfully terminate the calibration. This can be derived by solving Equation \eqref{eqn:P3_dpiC} for $\sRef$, and substituting \eqref{eqn:P3_dpiC} and \eqref{eqn:P3_lambda0} into \eqref{eqn:P3_dqC}; \citet{Liu2006} also follow this approach. 

Calibrating the model in this way may not be possible due to the following reasons:
In contrast to the other parameters and variables, the required reference marginal supply cost $\phiRef$ is difficult to obtain from data for a network of markets due to its dependency on $\qRef$ and the cost structure of the underlying network; Furthermore, $\MP_f(\lambdaRef$, $\qRef_f$, $\sRef$, $\phiRef_f$, $\etaRef$) is neither guaranteed to be in the interval $[0,1]$ for all $f \in \F$, nor properly defined for those traders whose sales $\qRef_f=0$. 

In the remainder of this section we derive three modules, based on which we then formulate a calibration algorithm overcoming these difficulties. 
\begin{itemize}
\item \mbox{Module I} (\mbox{Section \ref{sec:P3_DetOfPhi}}) introduces how the marginal costs $\phiRef$ of supplying $\qRef$ can be determined using the model.
\item \mbox{Module II} (\mbox{Section \ref{sec:P3_bounds}}) derives ranges on the anchor price $\lambdaO \in \lambdaORange$ and the price elasticity $\eta \in \etaRange$ for which the market power parameters $\MP_f(\lambdaO$, $\qRef_f$, $\sRef$, $\phiRef_f$, $\eta$) are in the interval $[0,1]$ for all $f \in \F$. 
\item If these ranges do not contain $\lambdaRef$ and $\etaRef$, we face a basic misalignment between reference data and market equilibria which the model can generate, and thus have to relax the constraints on the reference values. In this case we allow any $\lambdaRef \in [\underline{\lambdaRef},\overline{\lambdaRef}]$ and $\etaRef \in \etaRefRange$, and any equilibrium ($\lambdaRef$, $\qRef$, $\sRef$, $\phiRef$) to terminate the calibration, where $\underline{\lambdaRef}$, $\overline{\lambdaRef}$, $\underline{\etaRef}$, $\overline{\etaRef}$ are user-defined bounds on the reference values, for instance derived from uncertainty inherent to the reference data\footnote{As reference data usually comes with some uncertainty, we would argue that this approach is legitimate. Alternatively, the cost structure of the network, or the model altogether could be changed. In this work, however, we do not investigate those options.}.
If the ranges $\lambdaORange$ and $\etaRange$ are empty, or do not overlap with $\lambdaRefRange$ and $\etaRefRange$, we also allow changes in $\qRef$. \mbox{Module III} (\mbox{Section \ref{sec:P3_AlterQf}}) describes how a new $\qRef$ can be found, for which the ranges $\lambdaORange$ and $\etaRange$ are shifted towards $\lambdaRefRange$ and $\etaRefRange$.
\end{itemize}
Module I-III are integrated to a calibration algorithm in \mbox{Section \ref{sec:P3_DefCalbAlg}}. Notation is introduced as we proceed, and is summarized in \mbox{\ref{app:P3_NotationCalibAlg}}.

\subsection{Module I: Determining the marginal costs \texorpdfstring{$\phiRef{}$}{} of supplying the reference sales of the traders \texorpdfstring{$\qRef$}{}} \label{sec:P3_DetOfPhi}
We introduce 
\begin{equation}
0 =\qRef_f - q_{f}  \perp \xi_f \,\text{(free)}  \quad \forall f \in \mathcal{F}, \label{eqn:P3_depsFC}\\
\end{equation}
which are additional \ac{KKT} conditions from adding the constraint $q_{f}=\qRef_f$ to the original optimization problem of each trader $f$.
 $\xi_f$ are the shadow prices associated with these constraints and therefore also appear in Equation \eqref{eqn:P3_dqC}, which now reads
\begin{equation} 
0 \leq\, -\lambda -\MP_f \frac{d\IDF(s)}{ds} q_{f} + \phi_{f} + \xi_f \perp q_{f} \,\geq 0 \quad \forall f \in \mathcal{F}. \label{eqn:P3_dqCnew}
\end{equation}
We solve the augmented model comprising Equations \eqref{eqn:P3_dpiC}, \eqref{eqn:P3_depsFC}, and \eqref{eqn:P3_dqCnew}, choose parameters $\sO=\sRef$, $\lambdaO=\lambdaRef$, $\eta=\etaRef$, and any $\MP \succeq \vec{0}$, where $\vec{0}$ is the zero vector of appropriate dimension and $\succeq$ is interpreted component-wise. The variable $\xi_f$ corresponds to a tax/subsidy imposed on the trader $f$ to supply the quantity $\qRef_f$, and compensates for the potential mismatch of supply and demand in the market caused by our (arbitrary) choice of $\MP$. In equilibrium, the vector $\phi$ reveals the marginal costs of supplying $\qRef$, and we obtain $\phiRef := \phiStar$. 

\subsection{Module II: Determining admissible ranges for the anchor price \texorpdfstring{$\lambdaO$}{} and the price elasticity \texorpdfstring{$\eta$}{}} \label{sec:P3_bounds}
We obtain all values for the anchor price $\lambdaO$ and the price elasticity $\eta$ for which the market power parameters $\MP_f(\lambdaO$, $\qRef_f$, $\sRef$, $\phiRef_f$, $\eta$) are in the interval $[0,1]$ for all $f \in \F$ by exploring the dynamics inherent to the model. Note that $\F=\FPlus \cup \FZero = \lbrace f \in \{\F | \qRef_f > 0\} \rbrace \cup \lbrace f \in \{\F | \qRef_f = 0\} \rbrace$.

For $f \in \FPlus$, and provided $\MP = \MP(\lambdaO$, $\qRef$, $\sRef$, $\phiRef$, $\eta$), Equation \eqref{eqn:P3_dqC} implies that the following equation holds in equilibrium.  
\begin{equation} 
\lambdaO -\MP_f \frac{\lambdaO}{\sRef \cdot (-\eta)} \qRef_{f} = \phiRef_{f} \quad \forall f \in \FPlus \label{eqn:P3_Equate_q>0}
\end{equation} 
Since $\MP_f \frac{\lambdaO}{\sRef \cdot (-\eta)} \qRef_{f}$ is non-negative for all valid choices of $\lambdaO$, $\eta$, and $\MP_f \in [0,1]$, we conclude that the wholesale price in equilibrium $\lambdaStar = \lambdaO$ is at least as large as the highest marginal supply cost of the supplying traders.
\begin{equation} 
\lambdaO \geq \max \limits_{f \in \FPlus} \phiRef_f \label{eqn:P3_lambdaCalibMin}
\end{equation}
Furthermore, $\lambdaO$ is at most equal to marginal cost of the supplying trader with the lowest cost plus its maximum market power markup. 
\begin{equation} 
\lambdaO \leq \phiRef_f + 1 \cdot \frac{\lambdaO}{\sRef \cdot (-\eta)} \qRef_{f} \quad \forall f \in \FPlus. \label{eqn:P3_lambdaEtaAdmMax}
\end{equation}
We can reformulate condition \eqref{eqn:P3_lambdaEtaAdmMax} to provide an upper bound on $(-\eta)$ instead of $\lambdaO$:
\begin{equation} 
(-\eta) \leq \min \limits_{f \in \FPlus} \frac{\lambdaO}{\lambdaO- \phiRef_f} \cdot \frac{\qRef_f}{\sRef}. \label{eqn:P3_MinusEtaCalibMax}
\end{equation}

For $f \in \FZero$ Equation \eqref{eqn:P3_dqC} reads 
\begin{equation} 
\lambdaO \leq \phiRef_{f}, \label{eqn:P3_dqCnewEquilibriumq=0}
\end{equation}
limiting the price in equilibrium $\lambdaO$ to the lowest marginal cost among the non-supplying traders with access to the market:
\begin{equation} 
\lambdaO \leq \min \limits_{f \in \FZero} \phiRef_f.  \label{eqn:P3_lambdaCalibMax}
\end{equation}

Summarizing, if we limit $\MP$ to the interval $[\vec{0},\vec{1}]$ and consider the properties \eqref{eqn:P3_lambda0}, $\lambdaO$ and $\eta$ have to be in the ranges
\begin{equation} 
\underline{\lambdaO} := \max \limits_{f \in \FPlus} \phiRef_f \leq \lambdaO \leq \min \limits_{f \in \FZero} \phiRef_f =: \overline{\lambdaO}, \label{eqn:P3_lambdaCalibAdm}
\end{equation}
\begin{equation} 
\underline{\MinusEtaBounds} := 0 < (-\eta) \leq \min \limits_{f \in \FPlus} \frac{\lambdaO}{\lambdaO- \phiRef_f} \cdot \frac{\qRef_f}{\sRef} =: \overline{\MinusEtaBounds}. \label{eqn:P3_etaCalibAdm}
\end{equation}
If the boundary values of the admissible ranges satisfy $\overline{\lambdaO} \geq \underline{\lambdaO}$ and $\overline{\MinusEtaBounds} > \underline{\MinusEtaBounds}$, then the equilibrium ($\lambdaStar$, $\qStar$, $\sStar$, $\phiStar$) $=$ ($\lambdaO$, $\qRef$, $\sRef$, $\phiRef$) exists, and can be achieved by choosing any $\lambdaO$ and $\eta$ in the ranges $\lambdaORange$, $\etaRange$, setting $\MP = \MP(\lambdaO$, $\qRef_f$, $\sRef$, $\phiRef_f$, $\eta$) for all $f \in \FPlus$, and setting $\MP_f$ to any value in the interval [0,1] for all $f \in \FZero$; note that in all cases $\MP_f \in [0,1]$ as desired. 

\subsection{Module III: Improving the anchor price and price elasticity ranges \texorpdfstring{$\lambdaORange$ and $\etaRange$}{} by altering the reference sales of the traders} \label{sec:P3_AlterQf}
If the anchor price and price elasticity ranges $\lambdaORange$ and $\etaRange$ do not intersect with $\lambdaRefRange$ and $\etaRefRange$, we know that we cannot achieve the equilibrium ($\lambdaRef$, $\qRef$, $\sRef$, $\phiRef$) without violating $\lambdaO \in \lambdaRefRange$, $\eta \in \etaRefRange$, or $\MP_f \in [0,1]$ for some $f \in \F$. We follow a two-step procedure to find a new equilibrium ($\lambdaRef$, $\qNew$, $\sRef$, $\phiNew$) for which none of the parameters $\lambdaO$, $\eta$, and $\MP$ violate their respective intervals.

In a first step, we pick a $\lambdaO \in \lambdaRefRange$ and $\eta \in \etaRefRange$, which are as close as possible to $\lambdaORange$ and $\etaRange$, respectively. Furthermore, we fix the anchor consumption $\sO$ at $\sRef$. We calculate $\MP_f(\lambdaO$, $\qRef_f$, $\sO$, $\phiRef_f$, $\eta$) $=:\MPRef_f$ for all $f \in \FPlus$; these values can be outside the interval $[0,1]$. Furthermore, we define $\MPLim_f := \min(\max(\MPRef_f,0),1)$ for all $f \in \FPlus$, and set $\MPLim_f = \sum \limits_{f^\prime \in \FPlus} \frac{\MPLim_{f^\prime} \cdot \qRef_{f^\prime}}{\qRef_{f^\prime}}$ for all $f \in \FZero$, as reference data does not provide any information on the behavior of traders $f$ with $\qRef_f=0$; other approximations could be used as well. From Equation \eqref{eqn:P3_dqC} we estimate the sales in the new equilibrium:
\begin{subequations} \label{eqn:P3_qOIterPlusOneEstimate}
\begin{align}
\qEst_{f} &:= \frac{\lambdaO - \phiRef_{f}}{\lambdaO} \frac{\sO \cdot (-\eta)}{\MPLim_f} \quad \forall f \in \{\F| \MPLim_f>0\}, \\
\qEst_{f} &:= \qRef - \left|\frac{\lambdaO - \phiRef_{f}}{\lambdaO}\right| \frac{\sO \cdot (-\eta)}{1} \quad \forall f \in \{\F| \MPLim_f = 0\}.
\end{align}
\end{subequations}
Note that $\qEst_f = \qRef_f$ if $\MPRef_f$ is in the interval $[0,1]$. 
 
In a second step, we want to find a new equilibrium in the direction of $\qEst$. To this end, we augment the original model by
\begin{subequations} \label{eqn:P3_AugmentedModel2}
\begin{alignat}{4}
0 &\leq&\, \max \left(\qRef_{f}, \qEst_{f}\right) - q_{f} & \perp \xiFH_f &&\,\geq 0  \quad \forall f, \label{eqn:P3_depshiFC}\\
0 &\leq&\, -\min \left(\qRef_{f}, \qEst_{f}\right) + q_{f} & \perp \xiFL_f &&\,\geq 0  \quad \forall f, \label{eqn:P3_depsloFC}\\
0 &=&		\sO - s 					& \perp \chi &&\, \text{(free)}.   \label{eqn:P3_dchiC}
\end{alignat}
\end{subequations}
These complementarity constraints are the \ac{KKT} conditions of additional constraints in the traders' optimization problems. $\xiFH_f$, $\xiFL_f$, and $\chi$ are the shadow prices associated with these constraints and therefore also appear in Equation \eqref{eqn:P3_dqC}, which now reads
\begin{equation} 
0 \leq\, -\lambda -\MPLim_f \frac{d\IDF(s)}{ds} q_{f} + \phi_{f} + \xiFH_f- \xiFL_f+\chi \perp q_{f} \,\geq 0 \quad \forall f \in \mathcal{F}. \label{eqn:P3_dqCnew2}
\end{equation}
By solving the augmented model comprising Equations \eqref{eqn:P3_dpiC}, \eqref{eqn:P3_AugmentedModel2}, and \eqref{eqn:P3_dqCnew2}, with parameters $\sO=\sRef$, $\lambdaO$, $\eta$, and $\MP=\MPLim$, we obtain a new equilibrium ($\lambdaStar=\lambdaO$,$\qStar=:\qNew$,$\sStar=\sO$,$\phiStar=:\phiNew$). Note that $\lambdaStar=\lambdaO$ follows from the definition of the inverse demand function \eqref{eqn:P3_lambda0} and \mbox{constraint \eqref{eqn:P3_dchiC}}.

\subsection{Definition of the calibration algorithm} \label{sec:P3_DefCalbAlg}
From the previously introduced modules we can formulate a calibration algorithm for \ac{CV} models.
\begin{itemize}
\item Step 1 (Initialization): Derive reference parameter values $\lambdaRef$, $\etaRef$, and $\qRef$, as well as upper and lower bounds on these values $\underline{\lambdaRef}$, $\overline{\lambdaRef}$, $\underline{\etaRef}$, and $\overline{\etaRef}$, for instance from historical data. Set $\sO \equiv \sRef = \sum \limits_{f \in \F} \qRef_{f}$, set the iteration counter $i=1$, and define $\qRefIterOne:=\qRef$. Apply \mbox{Module I} (\mbox{Section \ref{sec:P3_DetOfPhi}}) to calculate $\phiRefIterOne$ from $\qRefIterOne$. 
\item Step 2: Calculate $\lambdaORangeIter$ and $\etaRangeIter$ from $\phiRefIter$ and $\qRefIter$ via Equations \eqref{eqn:P3_lambdaCalibAdm} and \eqref{eqn:P3_etaCalibAdm} (\mbox{Module II}).
\item Step 3 (Termination): If $\lambdaORangeIter$ and $\etaRangeIter$ are non-empty and intersect with $\lambdaRefRange$ and $\etaRefRange$:
Choose $\lambdaOIter \in \lambdaORangeIter \cap \lambdaRefRange$, $\etaIter \in \etaRangeIter \cap \etaRefRange$, calculate $\MPIter = \MP(\lambdaOIter,\qRefIter,\sRef,\phiRefIter,\etaIter)$, and terminate the calibration. 
Otherwise, move to Step 4.
\item Step 4 (Update): Apply \mbox{Module III} to obtain $\qRefIterPlusOne := \qNew$ and $\phiRefIterPlusOne$ $:=$ $\phiNew$. Update $i = i+1$, and move to \mbox{Step 2}.
\end{itemize}

Note that this algorithm is only guaranteed to terminate if the interval $\lambdaRefRange$ is chosen large enough with respect to all the reference values given. Unfortunately, ``large enough'' is difficult to quantify beforehand. As a consequence, the modeler might prefer to start with $\underline{\lambdaRef} = \overline{\lambdaRef} = \lambdaRef$, and gradually widen $\lambdaRefRange$ for those nodes and time periods stalling the algorithm. Furthermore, a less strict termination criterion can reduce the number of iterations; for instance $|\sStar-\sRef| \leq \TOL$ is suitable for most practical purposes, where $\sStar$ is the solution to the model comprising Equations \eqref{eqn:P3_dqC} and \eqref{eqn:P3_dpiC}, and $\TOL$ the maximum allowed deviation to the reference value.

\section{Numerical example} \label{sec:P3_NumericalExample}
We demonstrate the functionality of the proposed algorithm by applying it to the gas market model introduced by \citet{Baltensperger2015}. The model represents the \ac{EU} markets and their main suppliers over 2 periods (summer and winter), consists of 43 nodes and 247 arcs, and is represented by 9432 complementarity conditions (including Equations \eqref{eqn:P3_dqC} and \eqref{eqn:P3_dpiC}). The model equations and an exemplary model with two interconnected nodes are shown in \mbox{\ref{app:P3_ModelEq_Pic_Notation}}. The model and the calibration algorithm were implemented in MATLAB and solved by CPLEX. It takes a quad-core \unit[3.4]{GHz} CPU \unit[5.7]{seconds} on the average to compute one iteration of the algorithm, whereas the main computational burden originates from solving the augmented version of the model in Module III.

Historical data ($\lambdaData$, $\qData$, $\sData$, $\etaData$) was obtained from various sources \citep{ENTSO-G2012c, EuropeanCommission2013a, EuropeanCommission2013, EuropeanCommission, InternationalEnergyAgencyIEA2012c, InternationalEnergyAgencyIEA2013, Iea2013, Union2014, Lise2008, UnitedNations}.
As is often the case in practice, the total consumption did not match the reported sales of the traders $\sData \neq \sum \limits_{f \in \F} \qData_{f}$. A consistent set ($\lambdaRef$, $\qRef$, $\sRef \equiv \sum \limits_{f \in \F} \qRef_{f}$, $\etaRef$) was obtained by solving the augmented version of the model comprising Equations \eqref{eqn:P3_dpiC}, \eqref{eqn:P3_AugmentedModel2}, and \eqref{eqn:P3_dqCnew2} before initializing the algorithm. We used $\lambdaO_{nt} = \lambdaData_{nt}$, $\eta_{nt} = \etaData_{nt}$, $\MP_{fnt} = 1$ for all $f \in \F, n \in \N, t \in \T$, where $\N$ is the set of all nodes, and $\T$ the set of all time periods.
Equations \eqref{eqn:P3_depshiFC} and \eqref{eqn:P3_depsloFC} were altered to constrain $q_{fnt}$ to the interval $[\qData_{fnt}, \infty)$, if $\sum \limits_{f \in \mathcal{F}} \qData_{fnt} \leq \sData_{nt}$ and $\sum \limits_{n \in \mathcal{N}} \qData_{fnt} \leq \frac{\pData_{ft}}{1-\hat{\LOSS}_{ft}}$, where $\pData_{ft}$ is the total quantity produced by trader $f$ in period $t$, and $\hat{\LOSS}_{ft}$ an estimate of the lost fraction of gas until it reaches the consumers.
Otherwise, the $\qData_{fnt}$'s were scaled down to fit these inequalities, because $\sData_{nt}$ and $\pData_{ft}$ are among the most accurate figures available for the gas market.
In other situations, different parameters might be prioritized. We obtain ($\lambdaRef$, $\qRef$, $\sRef$, $\phiRef$) $=$ ($\lambdaStar=\lambdaData$, $\qStar$, $\sStar=\sData$, $\phiStar$).
For better readability, we drop subscripts $n$ and $t$ in the following.

In our example, we set $\underline{\lambdaRef} = \overline{\lambdaRef} = \lambdaRef$, $\underline{\etaRef} = \max(-\etaRef-0.2, 0.3)$, and $\overline{-\etaRef} = \min(-\etaRef+0.2, 1)$. If $\lambdaOIter$ was set to $\underline{\lambdaRef}$ in iteration $i$, we decreased $\underline{\lambdaRef}$ by $0.02 \cdot \lambdaRef$ in iteration $i+1$; we proceeded similarly when $\lambdaOIter$ hit the upper bound $\overline{\lambdaRef}$. The algorithm was set to terminate if neither $\underline{\lambdaRef}$ nor $\overline{\lambdaRef}$ were changed in an iteration, and $|\sStar-\sRef| \leq 0.5$ million cubic meters per day (\unitfrac[]{mcm}{d}). 

For the chosen bounds, the algorithm terminates in the \nth{10} iteration. The main characteristics of the solution are displayed in \mbox{Table \ref{tab:P3_calibration}}. Further details on the results are presented and discussed in \mbox{\ref{app:P3_AddResults}}. The top half of \mbox{Table \ref{tab:P3_calibration}} illustrates the deviations of the calibrated anchor values $\lambdaOIterEnd$ and $\etaIterEnd$ to their reference values. These deviations increase with increasing misalignment between reference data and market equilibria the model can generate. Furthermore, $|\qRefIterEnd-\qRef|$ increases with tighter bounds $\lambdaRefRange$ and $\etaRefRange$, and vice versa. Consequently, the modeler can distribute the deviations to the parameters of his choice by setting $\lambdaRefRange$ and $\etaRefRange$ accordingly.
Note that the infinite relative deviations in $\qRefIterEnd_f$ and $q_f$ originate from some $\qRefIterEnd_f>0$ and $q_f>0$ while $\qRef_f=0$. The low mean and median deviations indicate that the algorithm does not bias the results, for instance, towards a higher or lower average price level. 

The lower half of \mbox{Table \ref{tab:P3_calibration}} illustrates that the calibrated parameters $\lambdaOIterEnd$, $\etaIterEnd$, $\MP(\lambdaOIterEnd$, $\qRefIterEnd$, $\sO$, $\phiRefIterEnd$, $\etaIterEnd$) (and $\sO=\sRef$) indeed generate an equilibrium very close to the desired values ($\lambdaOIterEnd$, $\qRefIterEnd$, $\sO$, $\phiRefIterEnd$), which proves that the proposed algorithm works as intended when applied in practice.

\begin{table}[htbp]
\caption{The upper rows display the deviations between parameters $\lambdaOIterEnd$, $\qRefIterEnd_f$, and $-\etaIterEnd$, and their reference values. The lower rows show the deviations between the equilibrium values of the calibrated model and the calibrated parameters. Quantities are given in million cubic meters per day (\unitfrac[]{mcm}{d}), prices in thousand Euros per million cubic meters (\unitfrac[]{k\euro{}}{mcm}), and price elasticities are unitless.}
\label{tab:P3_calibration} 
\centering \small \renewcommand\arraystretch{1.5}
\addtolength{\tabcolsep}{-2pt}
\begin{tabu} to 0.9\linewidth{ccrclrclll}
\toprule
\multicolumn{2}{c}{\textbf{Parameter /} }				& \multicolumn{8}{c}{\textbf{Deviation between values}} \\
\multicolumn{2}{c}{\textbf{variable} }						& \textbf{min} \hspace*{-5pt}	&\hspace*{-5pt}\textbar\hspace*{-5pt}& \hspace*{-5pt} \textbf{max, abs.} 											& \textbf{min} \hspace*{-5pt}	&\hspace*{-5pt}\textbar\hspace*{-5pt}& \hspace*{-5pt} \textbf{max, rel.} 						& \textbf{mean} 															& \textbf{median} \\ 
\midrule
$\lambdaOIterEnd$ &	$\lambdaRef$ 				& -23.2 \hspace*{-5pt}&\hspace*{-5pt}\textbar\hspace*{-5pt}& \hspace*{-5pt} \unitfrac[39.4]{k\euro{}}{mcm} 	& -7.56 \hspace*{-5pt}&\hspace*{-5pt}\textbar\hspace*{-5pt}& \hspace*{-5pt} \unit[16.8]{\%} 			& \unitfrac[-0.34]{k\euro{}}{mcm} 	& \unitfrac[0.00]{k\euro{}}{mcm}\\
$-\etaIterEnd$		& $-\etaRef$					& -0.20	\hspace*{-5pt}&\hspace*{-5pt}\textbar\hspace*{-5pt}& \hspace*{-5pt} \unit[0.20]{} 									& -37.4 \hspace*{-5pt}&\hspace*{-5pt}\textbar\hspace*{-5pt}& \hspace*{-5pt} \unit[59.0]{\%} 			& \unit[-0.01]{} 										& \unit[0.00]{}\\
$\qRefIterEnd_{f}$& $\qRef_{f}$	 				& -17.5	\hspace*{-5pt}&\hspace*{-5pt}\textbar\hspace*{-5pt}& \hspace*{-5pt} \unitfrac[17.5]{mcm}{d} 				& -100 	\hspace*{-5pt}&\hspace*{-5pt}\textbar\hspace*{-5pt}& \hspace*{-5pt} \unit[$\infty$]{\%} 	& \unitfrac[0.00]{mcm}{d} 					& \unitfrac[0.00]{mcm}{d}\\
\midrule
$\sStar$ 							& $\sRef$ 						&	-0.03 \hspace*{-5pt}&\hspace*{-5pt}\textbar\hspace*{-5pt}& \hspace*{-5pt} \unitfrac[0.41]{mcm}{d} 				& -0.03 \hspace*{-5pt}&\hspace*{-5pt}\textbar\hspace*{-5pt}& \hspace*{-5pt} \unit[0.25]{\%} 			& \unitfrac[0.00]{mcm}{d}						& \unitfrac[0.00]{mcm}{d}\\
$\lambdaStar$ 				& $\lambdaOIterEnd$		&	-3.15 \hspace*{-5pt}&\hspace*{-5pt}\textbar\hspace*{-5pt}& \hspace*{-5pt} \unitfrac[0.26]{k\euro{}}{mcm} 	& -0.83 \hspace*{-5pt}&\hspace*{-5pt}\textbar\hspace*{-5pt}& \hspace*{-5pt} \unit[0.07]{\%} 			& \unitfrac[-0.05]{k\euro{}}{mcm} 	& \unitfrac[0.05]{k\euro{}}{mcm}\\
$\qStar_f$ 						& $\qRefIterEnd_{f}$	& -0.07	\hspace*{-5pt}&\hspace*{-5pt}\textbar\hspace*{-5pt}& \hspace*{-5pt} \unitfrac[0.05]{mcm}{d} 				&-39.3 	\hspace*{-5pt}&\hspace*{-5pt}\textbar\hspace*{-5pt}& \hspace*{-5pt} \unit[$\infty$]{\%} 	& \unitfrac[0.00]{mcm}{d} 					& \unitfrac[0.00]{mcm}{d}\\
\bottomrule 
\end{tabu}
\end{table}

Note that our algorithm reproduces the rudimentary algorithm introduced in Section \ref{sec:P3_ProblemStatement} in the first iteration, if $\lambdaRefRange := [0,\infty)$ and $\etaRefRange := (-\infty,0)$ are chosen. In this light, our algorithm is a significant extension, as it finds a solution with $\MP$ in the interval $[\vec{0},\vec{1}]$, and highlights the interplay between the parameters, which allows the user to take informed decisions on how the mismatch between model and reality is resolved.

\section{Summary and Outlook} \label{sec:P3_Conclusions}
We propose an iterative algorithm to calibrate conjectural variations models for a network of markets. The algorithm builds upon three modules. In \mbox{Module I}, the marginal supply costs $\phi_f$ are derived for all traders $f$ from the sales $\qRef$ of the traders and the cost structure of the network. In \mbox{Module II}, all anchor prices $\lambdaO \in \lambdaORange$ and price elasticities $\eta \in \etaRange$, which enable the equilibrium of the model to be in a certain point, are explored. In \mbox{Module III}, a new consistent and physically feasible set of sales $\qNew$ is obtained, which changes the marginal costs of supply $\phi$ such that the ranges $\lambdaORange$ and $\etaRange$ are shifted in a desired direction. 

The example indicates that the algorithm is able to calibrate a network of markets to real world data fast after a small number of iterations. The presented algorithm is broadly applicable, since the assumptions on which it is based are very mild and common in spatial partial equilibrium modeling. Furthermore, the algorithm can easily be configured to specific purposes, since the reference values and their tolerances can be chosen by the user. Finally, the algorithm can easily be expanded; for instance, one could follow the approach taken by \citet{Huppmann2014} and (manually) adjust the underlying marginal cost parameters of specific system services to provoke a more favorable outcome of the sales per trader.

Future work should emphasize on extending the update step of the algorithm, such that the \textit{sizes} of the ranges of the various reference values can be weighed against each other and adjusted accordingly. In order to further increase the plausibility of the outcome with respect to economic theory, we also aim to improve how the algorithm sets these ranges relatively to each other. For instance, the market power parameters could be restricted to be larger in the high demand season than in the low demand season. Furthermore, we expect to achieve a reduction of iterations by exploiting the problem structure even further, particularly when determining the new set of sales per trader $\qNew$ in the update step.

\section*{Acknowledgements}
\addcontentsline{toc}{section}{Acknowledgements}
The authors are grateful to Prof.~Ruud Egging, for his insightful comments on earlier versions of the manuscript.

\clearpage
\appendix
\section{Calibration algorithm specific notation} \label{app:P3_NotationCalibAlg}

\begin{table}[htbp]
\caption{Variables, parameters, and functions.}
\label{tab:P3_vars}
\centering \small \renewcommand\arraystretch{1.5}
\begin{tabu} to 0.9\linewidth{lX}
\toprule 
\textbf{Abbr.} & \textbf{Explanation} \\
\midrule 
$q_f$& Gas shipment of \mbox{trader $f$}\\
$p_f$& Total production of \mbox{trader $f$}\\
$s$& Total consumption\\
$\eta$& Price elasticity of demand\\
$\MP_f$& Market power parameter of \mbox{trader $f$}\\
$\lambda$& Wholesale market price\\
$\IDF$& Inverse demand function \\
$\phi_f$& Marginal cost of trader $f$ to supply \mbox{quantity $q_f$}\\
$\xiFH_f$ & Shadow price of enforcing an upper level on trader $f$'s \mbox{sales $q_f$}\\
$\xiFL_f$ & Shadow price of enforcing an lower level on on trader $f$'s \mbox{sales $q_f$}\\
$\chi$ & Shadow price of enforcing consumption \mbox{level $\sO$}\\ 
\bottomrule
\end{tabu}
\end{table}

\begin{table}[htbp]
\caption{Sub- and superscripts.}
\label{tab:P3_subsup}
\centering \small \renewcommand\arraystretch{1.5}
\begin{tabu} to 0.9\linewidth{lX}
\toprule 
\textbf{Abbr.} & \textbf{Explanation} \\
\midrule 
$(\cdot)_f$& Parameter/variable of \mbox{trader $f$}\\
$(\cdot)_n$& Parameter/variable in \mbox{node $n$}\\
$(\cdot)_t$& Parameter/variable in time \mbox{period $t$}\\ 
\midrule
$(\cdot)^0$ & Anchor value for inverse demand function \\
$(\cdot)^{\ast}$ & Value of variable in equilibrium \\
$(\cdot)^\Ref$& Reference value for calibration\\
$(\cdot)^\Data$& Reference data, for instance based on historical values\\
$(\cdot)^\Est$& Estimated variable\\
$(\cdot)^\new$& Updated variable\\
$(\cdot)^i$ & Parameter/variable in \mbox{iteration $i$} \\
\midrule
$\overline{(\cdot)}$ & upper bound on parameter/variable\\
$\underline{(\cdot)}$ & lower bound on parameter/variable\\
\bottomrule
\end{tabu}
\end{table}

\begin{table}[htbp]
\caption{Models and sets.}
\label{tab:P3_models}
\centering \small \renewcommand\arraystretch{1.5}
\begin{tabu} to 0.9\linewidth{lX}
\toprule 
\textbf{Abbr.} & \textbf{Explanation} \\
\midrule 
$\F$ & Set of all traders \\ 
$\FPlus$ & Set of traders selling gas; $\FPlus = \{\F | q_f > 0\}$\\
$\FZero$ & Set of traders present in the market but not selling gas; $\FZero = \{\F | q_f = 0\}$\\
$\N$ & Set of all nodes/markets \\ 
$\T$ & Set of all time periods \\ 
\bottomrule
\end{tabu}
\end{table}

\clearpage
\section{In-depth analysis of results} \label{app:P3_AddResults}

As \mbox{Figure \ref{fig:P3_qOIterDev}} depicts, the supplier's market shares after calibration $\frac{\qIterEnd_f}{\sRef}$ (lowermost bar in each country) are close to the reference values $\frac{\qRefIterOne_f}{\sRef}$ (topmost bar in each country) for most countries. We find the largest deviations between $\frac{\qIterEnd_f}{\sRef}$ and $\frac{\qRefIterOne_f}{\sRef}$ in countries with rather low overall consumption such as \BG{}, the \CZ{}, \PT{}, \SK{}, \SI{}, \SE{}, and \CH{}; see \mbox{Table \ref{tab:P3_OtherCalibValues}} for consumption values $\sRef$. At the same time, deviations are comparably low for the largest \ac{EU} consumers \IT{}, \DE{}, \tUK{}, \FR{}, and \tNL{}. We explain the dependency on the market size as follows: The calibration aims at finding a combination of $\phiRef$, $\lambdaO$, and $\qRef$, which produces an equilibrium in itself. Therefore, the algorithm shifts the market shares $\frac{\qRef_f}{\sRef}$ for some $f$ over the course of the iterations (Equation \eqref{eqn:P3_qOIterPlusOneEstimate}). These impact the marginal supply costs $\phiRef$ of all traders, but first and foremost of those traders $f$ whose market shares $\frac{\qRef_f}{\sRef}$ are shifted. Once $\phiRef_f$ leads to intersecting ranges $\lambdaORange$ and $\lambdaRefRange$, Equation \eqref{eqn:P3_qOIterPlusOneEstimate} outputs $\qEst = \qRef$, and $\frac{\qRef_f}{\sRef}$ stops shifting. Consequently, the largest shifts in market shares occur in those countries in which changes in $\frac{\qRef_f}{\sRef}$ have the lowest impact on the marginal costs $\phiRef$ and therefore the range $\lambdaORange$, which are countries with low consumption compared to the total production of their supplier, the available pipeline capacity for imports, and the regasification capacity. 

\begin{figure*}[htbp]  
\centering
\includegraphics[width=0.9\textwidth]{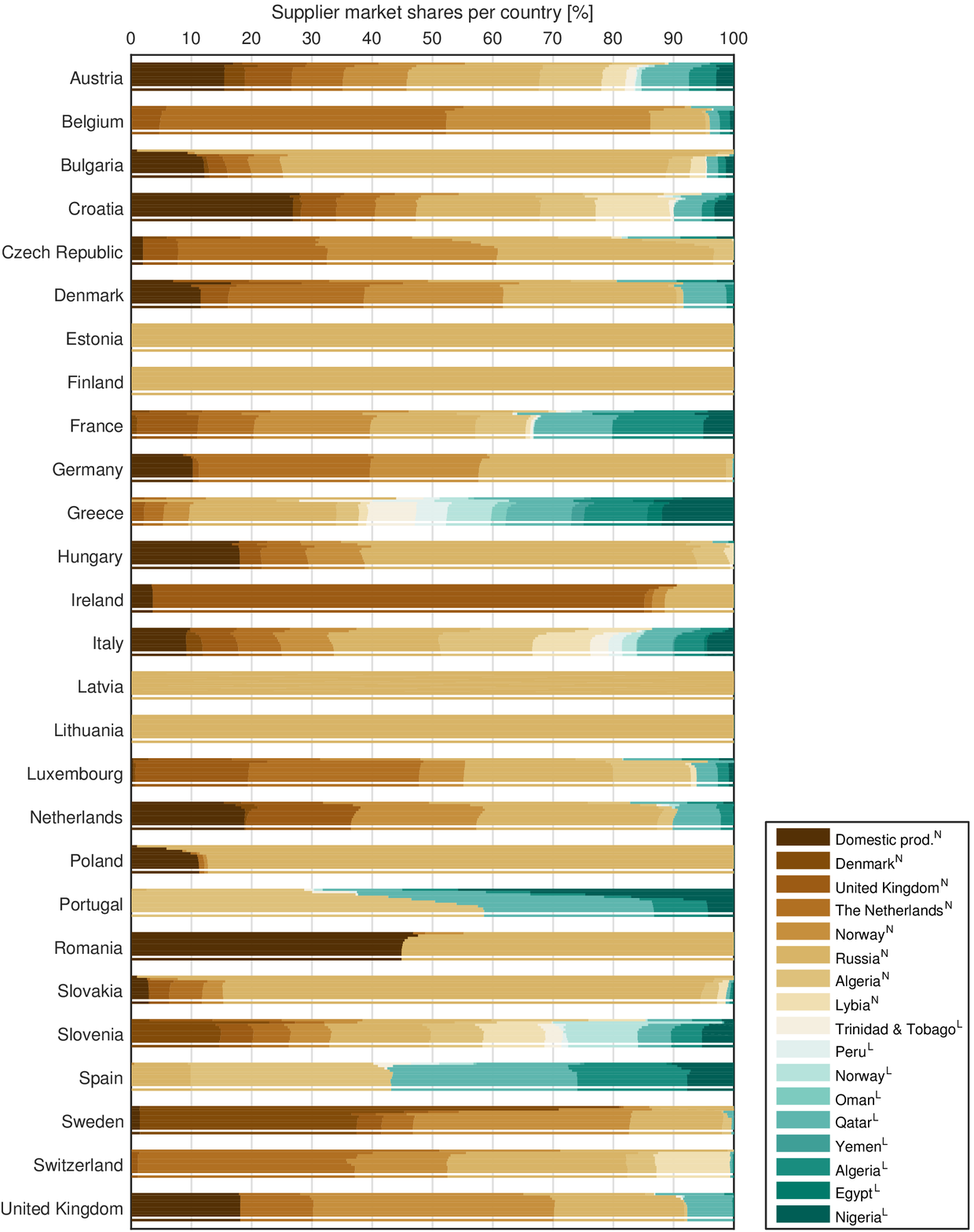}
\end{figure*}

\begin{figure*}[htbp]
\caption{Development of the market shares of all traders in all markets. For each country, a group of 10+1 bars is displayed, and each bar is divided into at most 17 sections. From top to bottom, the bars show the market shares $\frac{\qRefIter_f}{\sRef}$ of the suppliers in the iterations $i \in \{1, \ldots, 10\}$. The lowermost bar represents the market shares in the calibrated equilibrium $\frac{\qIterEnd_f}{\sRef}$.}
\label{fig:P3_qOIterDev} 
\end{figure*}

\begin{table}[htbp]
\caption{Reference consumption $\sRef$, reference price $\lambdaRef$, and reference price elasticity $\etaRef$, as well as the calibrated anchor consumption $\sO$, price $\lambdaOIterEnd$ and price elasticity $\etaIterEnd$. Quantities are given in million cubic meters per day (\unitfrac[]{mcm}{d}), prices in thousand Euros per million cubic meters (\unitfrac[]{k\euro{}}{mcm}), and price elasticities are unitless. $\lambdaOIterEnd$ and $\etaIterEnd$ above their respective references are colored red with increasing absolute value, and below below their respective references blue with decreasing value, to visually underpin our findings.}
\label{tab:P3_OtherCalibValues}
\centering \small \renewcommand\arraystretch{1.5}
\addtolength{\tabcolsep}{-4pt}
\begin{scriptsize}
\begin{tabu} to 0.9\linewidth{lcccccccccc}
		\toprule
          & \multicolumn{5}{c}{\textbf{October-March}} & \multicolumn{5}{c}{\textbf{April-September}} \\
          & \multicolumn{1}{c}{$\sRef$, $\sO$} & \multicolumn{1}{c}{$\lambdaRef$} & \multicolumn{1}{c}{$\lambdaOIterEnd$}  & \multicolumn{1}{c}{$\etaRef$} & \multicolumn{1}{c}{$\etaIterEnd$} & \multicolumn{1}{c}{$\sRef$, $\sO$} & \multicolumn{1}{c}{$\lambdaRef$} & \multicolumn{1}{c}{$\lambdaOIterEnd$}  & \multicolumn{1}{c}{$\etaRef$} & \multicolumn{1}{c}{$\etaIterEnd$} \\
					& \multicolumn{1}{c}{[\unitfrac[]{mcm}{d}]} & \multicolumn{1}{c}{[\unitfrac[]{k\euro{}}{mcm}]} & \multicolumn{1}{c}{[\unitfrac[]{k\euro{}}{mcm}]} & \multicolumn{1}{c}{[-]} & \multicolumn{1}{c}{[-]} & \multicolumn{1}{c}{[\unitfrac[]{mcm}{d}]} & \multicolumn{1}{c}{[\unitfrac[]{k\euro{}}{mcm}]} & \multicolumn{1}{c}{[\unitfrac[]{k\euro{}}{mcm}]} & \multicolumn{1}{c}{[-]} & \multicolumn{1}{c}{[-]} \\
    \midrule
    \AT{} & 34    & 354   & \ValuesColoredBR{90}349 & -0.47 & \ValuesColoredBR{57}-0.35 & 16    & 354   & \ValuesColoredBR{100}354 & -0.47 & \ValuesColoredBR{39}-0.30 \\
    \BE{} & 62    & 315   & \ValuesColoredBR{74}302 & -0.44 & \ValuesColoredBR{62}-0.34 & 39    & 315   & \ValuesColoredBR{99}315 & -0.44 & \ValuesColoredBR{47}-0.30 \\
    \BG{} & 10    & 413   & \ValuesColoredBR{71}392 & -0.43 & \ValuesColoredBR{48}-0.30 & 6     & 413   & \ValuesColoredBR{66}389 & -0.43 & \ValuesColoredBR{50}-0.30 \\
    \HR{} & 11    & 376   & \ValuesColoredBR{107}380 & -0.44 & \ValuesColoredBR{100}-0.44 & 7     & 376   & \ValuesColoredBR{95}372 & -0.44 & \ValuesColoredBR{90}-0.41 \\
		\CZ{} & 34    & 291   & \ValuesColoredBR{145}313 & -0.34 & \ValuesColoredBR{100}-0.34 & 13    & 291   & \ValuesColoredBR{100}291 & -0.34 & \ValuesColoredBR{200}-0.54 \\
    \DK{} & 13    & 305   & \ValuesColoredBR{104}307 & -0.40 & \ValuesColoredBR{100}-0.40 & 8     & 305   & \ValuesColoredBR{100}305 & -0.40 & \ValuesColoredBR{184}-0.60 \\
    \EE{} & 3     & 381   & \ValuesColoredBR{100}381 & -0.43 & \ValuesColoredBR{181}-0.63 & 2     & 381   & \ValuesColoredBR{100}381 & -0.43 & \ValuesColoredBR{178}-0.63 \\
    \FI{} & 14    & 381   & \ValuesColoredBR{89}373 & -0.50 & \ValuesColoredBR{100}-0.50 & 9     & 381   & \ValuesColoredBR{89}373 & -0.50 & \ValuesColoredBR{100}-0.50 \\
    \FR{} & 189   & 333   & \ValuesColoredBR{100}333 & -0.32 & \ValuesColoredBR{100}-0.32 & 71    & 333   & \ValuesColoredBR{96}331 & -0.32 & \ValuesColoredBR{93}-0.30 \\
    \DE{} & 287   & 305   & \ValuesColoredBR{100}305 & -0.41 & \ValuesColoredBR{185}-0.61 & 161   & 305   & \ValuesColoredBR{93}301 & -0.41 & \ValuesColoredBR{54}-0.30 \\
		\GR{} & 14    & 383   & \ValuesColoredBR{98}382 & -0.60 & \ValuesColoredBR{46}-0.42 & 12    & 383   & \ValuesColoredBR{72}365 & -0.60 & \ValuesColoredBR{63}-0.47 \\
    \HU{} & 42    & 337   & \ValuesColoredBR{111}343 & -0.35 & \ValuesColoredBR{100}-0.35 & 17    & 337   & \ValuesColoredBR{83}327 & -0.35 & \ValuesColoredBR{77}-0.30 \\
    \IE{} & 15    & 298   & \ValuesColoredBR{117}306 & -0.57 & \ValuesColoredBR{103}-0.58 & 13    & 298   & \ValuesColoredBR{94}295 & -0.57 & \ValuesColoredBR{42}-0.38 \\
    \IT{} & 286   & 376   & \ValuesColoredBR{99}375 & -0.47 & \ValuesColoredBR{39}-0.31 & 161   & 376   & \ValuesColoredBR{100}376 & -0.47 & \ValuesColoredBR{39}-0.30 \\
    \LV{} & 6     & 377   & \ValuesColoredBR{100}377 & -0.35 & \ValuesColoredBR{200}-0.55 & 2     & 377   & \ValuesColoredBR{100}377 & -0.35 & \ValuesColoredBR{197}-0.55 \\
		\LT{} & 9     & 421   & \ValuesColoredBR{100}421 & -0.53 & \ValuesColoredBR{165}-0.73 & 5     & 421   & \ValuesColoredBR{100}421 & -0.53 & \ValuesColoredBR{163}-0.73 \\
    \LU{} & 4     & 315   & \ValuesColoredBR{97}314 & -0.45 & \ValuesColoredBR{42}-0.30 & 3     & 315   & \ValuesColoredBR{97}314 & -0.45 & \ValuesColoredBR{100}-0.45 \\
    \TNL{} & 142   & 301   & \ValuesColoredBR{134}318 & -0.46 & \ValuesColoredBR{100}-0.46 & 80    & 301   & \ValuesColoredBR{109}306 & -0.46 & \ValuesColoredBR{100}-0.46 \\
    \PL{} & 55    & 301   & \ValuesColoredBR{55}278 & -0.37 & \ValuesColoredBR{95}-0.35 & 34    & 301   & \ValuesColoredBR{66}284 & -0.37 & \ValuesColoredBR{70}-0.30 \\
    \PT{} & 14    & 275   & \ValuesColoredBR{185}315 & -0.46 & \ValuesColoredBR{100}-0.46 & 14    & 275   & \ValuesColoredBR{158}302 & -0.46 & \ValuesColoredBR{100}-0.46 \\
    \RO{} & 42    & 177   & \ValuesColoredBR{200}207 & -0.43 & \ValuesColoredBR{100}-0.43 & 31    & 177   & \ValuesColoredBR{100}177 & -0.43 & \ValuesColoredBR{179}-0.63 \\
    \SK{} & 21    & 345   & \ValuesColoredBR{80}333 & -0.42 & \ValuesColoredBR{50}-0.30 & 8     & 345   & \ValuesColoredBR{79}333 & -0.42 & \ValuesColoredBR{51}-0.30 \\
		\SI{} & 3     & 376   & \ValuesColoredBR{87}368 & -0.40 & \ValuesColoredBR{91}-0.38 & 2     & 376   & \ValuesColoredBR{66}354 & -0.40 & \ValuesColoredBR{67}-0.32 \\
    \ES{} & 102   & 294   & \ValuesColoredBR{144}316 & -0.45 & \ValuesColoredBR{100}-0.45 & 85    & 294   & \ValuesColoredBR{117}303 & -0.45 & \ValuesColoredBR{100}-0.45 \\
		\SE{} & 5     & 305   & \ValuesColoredBR{108}309 & -0.53 & \ValuesColoredBR{100}-0.53 & 2     & 305   & \ValuesColoredBR{100}305 & -0.53 & \ValuesColoredBR{37}-0.33 \\
		\CH{} & 11    & 305   & \ValuesColoredBR{100}305 & -0.30 & \ValuesColoredBR{98}-0.30 & 4     & 305   & \ValuesColoredBR{97}303 & -0.30 & \ValuesColoredBR{100}-0.30 \\
		\UK{} & 256   & 298   & \ValuesColoredBR{107}301 & -0.44 & \ValuesColoredBR{100}-0.44 & 158   & 298   & \ValuesColoredBR{102}299 & -0.44 & \ValuesColoredBR{100}-0.44 \\
    \bottomrule
\end{tabu}%
\end{scriptsize}
\end{table}

Furthermore, we observe that the changes from $\frac{\qRefIterPlusOne_f}{\sRef}$ to $\frac{\qRefIter_f}{\sRef}$ are generally the largest in the first few iterations and decrease thereafter. Large changes indicate that the anchor price $\lambdaOIter$ changes, since this impacts the equilibrium heavily. This can be confirmed with help of \mbox{Table \ref{tab:P3_OtherCalibValues}}: for the example of \PT{} in October-March, $\frac{\lambdaOIterEnd}{\lambdaRef} = \frac{\unitfrac[315]{k\euro{}}{mcm}}{\unitfrac[275]{k\euro{}}{mcm}} \approx 1.15$, and as we increase $\overline{\lambdaRef}$ by $\unit[2]{\%}$ in the iteration $i+1$ if $\lambdaOIter = \overline{\lambdaRef}$ (\mbox{Section \ref{sec:P3_NumericalExample}}), it is safe to assume that $\lambdaOIterPlusOne = 1 + 0.02 \cdot i \cdot \lambdaRef$ for iterations $i \in \{1, \ldots, 7\}$, and only stagnates thereafter. This is confirmed by the stagnating shares of \PT{} after iteration $i=7$ (\mbox{Figure \ref{fig:P3_qOIterDev}}). 

We learn from \mbox{Table \ref{tab:P3_OtherCalibValues}} that anchor prices $\lambdaOIterEnd$ are overall slightly higher in the European winter than during summer. This intuitively makes sense, since consumption $\sRef$ is clearly higher in winter than in summer, while production capacities remain unchanged, and therefore gas is scarcer. For the price elasticities $\etaIterEnd$ we do not observe such a clear trend, but note that the $\etaIterEnd$ remain close to reference values $\etaRef$ for the larger consumers \UK{}, \FR{}, \tNL{}, and \ES{}. In \IT{}, $(-\etaIterEnd)$ is very low in both seasons. This is a consequence of the high wholesale prices $\lambdaRef$ compared to the marginal supply costs $\phiRef$ of the traders (Equation \eqref{eqn:P3_MinusEtaCalibMax}). The algorithm concludes from this difference that consumers in \IT{} are willing to buy gas at high prices and therefore assigns low $(-\eta)$ and high market power parameters (\mbox{Tables \ref{tab:P3_CVvaluesOC}} \mbox{and \ref{tab:P3_CVvaluesAP}}), which coincides with the interpretation an economist would make when analyzing such a situation. In \DE{}, we face a special situation: from October to March, $(-\eta)$ is capped at its maximum value, while the opposite is true for April to September. As previously, this result originates from the given reference data. We do not judge at this point whether the outcome is realistic from an economic point of view; instead we emphasize that the algorithm allows the modeler to set the ranges $\etaRefRange$ such that the outcome is defensible.

\begin{table}[htbp]
\caption{Market power parameter values for all traders (T) in all markets (C) in October-March. The country abbreviations are given in \mbox{Table \ref{tab:P3_CountryAbbr}}. \DZ{} and \NO{} supply the \ac{EU} via the pipeline network and via \ac{LNG} shipments. This is distinguished by ${}^N$ and ${}^L$, respectively. Figures above 0.5 are colored red with increasing value, and below 0.5 green with decreasing value, to visually underpin our findings.}
\label{tab:P3_CVvaluesOC}
\centering \small \renewcommand\arraystretch{1.5}
\addtolength{\tabcolsep}{-4pt}
\begin{scriptsize}
\begin{tabu} to 0.9\linewidth{crrrrrrrrrrrrrrrrr}
    \toprule
\backslashbox{C}{T} & AZ    & DK    & DZ$^L$    & DZ$^N$   & EG    & GB    & LY    & NG    & NL    & NO$^L$    & NO$^N$   & OM    & PE    & QA    & RU    & TT    & YE \\
    \midrule
    AT    & \ValuesColored{92}0.46 & \ValuesColored{82}0.41 & \ValuesColored{80}0.40 & \ValuesColored{84}0.42 & \ValuesColored{92}0.46 & \ValuesColored{76}0.38 & \ValuesColored{200}1 & \ValuesColored{166}0.83 & \ValuesColored{86}0.43 & \ValuesColored{0}0 & \ValuesColored{82}0.41 & \ValuesColored{92}0.46 & \ValuesColored{92}0.46 & \ValuesColored{34}0.17 & \ValuesColored{118}0.59 & \ValuesColored{44}0.22 & \ValuesColored{92}0.46 \\
    BE    & \ValuesColored{18}0.09 & \ValuesColored{18}0.09 & \ValuesColored{18}0.09 & \ValuesColored{18}0.09 & \ValuesColored{18}0.09 & \ValuesColored{0}0 & \ValuesColored{18}0.09 & \ValuesColored{18}0.09 & \ValuesColored{2}0.01 & \ValuesColored{18}0.09 & \ValuesColored{6}0.03 & \ValuesColored{18}0.09 & \ValuesColored{18}0.09 & \ValuesColored{18}0.09 & \ValuesColored{200}1 & \ValuesColored{18}0.09 & \ValuesColored{18}0.09 \\
    BG    & \ValuesColored{82}0.41 & \ValuesColored{200}1 & \ValuesColored{200}1 & \ValuesColored{200}1 & \ValuesColored{82}0.41 & \ValuesColored{200}1 & \ValuesColored{200}1 & \ValuesColored{168}0.84 & \ValuesColored{176}0.88 & \ValuesColored{82}0.41 & \ValuesColored{132}0.66 & \ValuesColored{82}0.41 & \ValuesColored{82}0.41 & \ValuesColored{156}0.78 & \ValuesColored{38}0.19 & \ValuesColored{82}0.41 & \ValuesColored{82}0.41 \\
    HR    & \ValuesColored{138}0.69 & \ValuesColored{200}1 & \ValuesColored{166}0.83 & \ValuesColored{164}0.82 & \ValuesColored{138}0.69 & \ValuesColored{110}0.55 & \ValuesColored{72}0.36 & \ValuesColored{70}0.35 & \ValuesColored{128}0.64 & \ValuesColored{138}0.69 & \ValuesColored{152}0.76 & \ValuesColored{138}0.69 & \ValuesColored{138}0.69 & \ValuesColored{50}0.25 & \ValuesColored{166}0.83 & \ValuesColored{138}0.69 & \ValuesColored{138}0.69 \\
    CZ    & \ValuesColored{28}0.14 & \ValuesColored{28}0.14 & \ValuesColored{28}0.14 & \ValuesColored{0}0 & \ValuesColored{28}0.14 & \ValuesColored{28}0.14 & \ValuesColored{28}0.14 & \ValuesColored{28}0.14 & \ValuesColored{6}0.03 & \ValuesColored{28}0.14 & \ValuesColored{10}0.05 & \ValuesColored{28}0.14 & \ValuesColored{28}0.14 & \ValuesColored{28}0.14 & \ValuesColored{56}0.28 & \ValuesColored{28}0.14 & \ValuesColored{28}0.14 \\
    DK    & \ValuesColored{28}0.14 & \ValuesColored{16}0.08 & \ValuesColored{28}0.14 & \ValuesColored{28}0.14 & \ValuesColored{28}0.14 & \ValuesColored{28}0.14 & \ValuesColored{28}0.14 & \ValuesColored{28}0.14 & \ValuesColored{0}0 & \ValuesColored{28}0.14 & \ValuesColored{16}0.08 & \ValuesColored{28}0.14 & \ValuesColored{28}0.14 & \ValuesColored{28}0.14 & \ValuesColored{64}0.32 & \ValuesColored{28}0.14 & \ValuesColored{28}0.14 \\
    EE    & \ValuesColored{62}0.31 & \ValuesColored{62}0.31 & \ValuesColored{62}0.31 & \ValuesColored{62}0.31 & \ValuesColored{62}0.31 & \ValuesColored{62}0.31 & \ValuesColored{62}0.31 & \ValuesColored{62}0.31 & \ValuesColored{62}0.31 & \ValuesColored{62}0.31 & \ValuesColored{62}0.31 & \ValuesColored{62}0.31 & \ValuesColored{62}0.31 & \ValuesColored{62}0.31 & \ValuesColored{62}0.31 & \ValuesColored{62}0.31 & \ValuesColored{62}0.31 \\
    FI    & \ValuesColored{48}0.24 & \ValuesColored{48}0.24 & \ValuesColored{48}0.24 & \ValuesColored{48}0.24 & \ValuesColored{48}0.24 & \ValuesColored{48}0.24 & \ValuesColored{48}0.24 & \ValuesColored{48}0.24 & \ValuesColored{48}0.24 & \ValuesColored{48}0.24 & \ValuesColored{48}0.24 & \ValuesColored{48}0.24 & \ValuesColored{48}0.24 & \ValuesColored{48}0.24 & \ValuesColored{48}0.24 & \ValuesColored{48}0.24 & \ValuesColored{48}0.24 \\
		FR    & \ValuesColored{54}0.27 & \ValuesColored{38}0.19 & \ValuesColored{22}0.11 & \ValuesColored{24}0.12 & \ValuesColored{54}0.27 & \ValuesColored{42}0.21 & \ValuesColored{54}0.27 & \ValuesColored{52}0.26 & \ValuesColored{56}0.28 & \ValuesColored{54}0.27 & \ValuesColored{80}0.40 & \ValuesColored{54}0.27 & \ValuesColored{54}0.27 & \ValuesColored{22}0.11 & \ValuesColored{108}0.54 & \ValuesColored{196}0.98 & \ValuesColored{54}0.27 \\
    DE    & \ValuesColored{60}0.30 & \ValuesColored{60}0.30 & \ValuesColored{60}0.30 & \ValuesColored{60}0.30 & \ValuesColored{60}0.30 & \ValuesColored{60}0.30 & \ValuesColored{60}0.30 & \ValuesColored{60}0.30 & \ValuesColored{6}0.03 & \ValuesColored{60}0.30 & \ValuesColored{30}0.15 & \ValuesColored{60}0.30 & \ValuesColored{60}0.30 & \ValuesColored{60}0.30 & \ValuesColored{92}0.46 & \ValuesColored{60}0.30 & \ValuesColored{60}0.30 \\
    GR    & \ValuesColored{124}0.62 & \ValuesColored{124}0.62 & \ValuesColored{122}0.61 & \ValuesColored{132}0.66 & \ValuesColored{200}1 & \ValuesColored{160}0.80 & \ValuesColored{124}0.62 & \ValuesColored{88}0.44 & \ValuesColored{136}0.68 & \ValuesColored{108}0.54 & \ValuesColored{130}0.65 & \ValuesColored{130}0.65 & \ValuesColored{142}0.71 & \ValuesColored{146}0.73 & \ValuesColored{124}0.62 & \ValuesColored{120}0.60 & \ValuesColored{200}1 \\
    HU    & \ValuesColored{60}0.30 & \ValuesColored{60}0.30 & \ValuesColored{60}0.30 & \ValuesColored{52}0.26 & \ValuesColored{60}0.30 & \ValuesColored{16}0.08 & \ValuesColored{0}0 & \ValuesColored{60}0.30 & \ValuesColored{22}0.11 & \ValuesColored{60}0.30 & \ValuesColored{32}0.16 & \ValuesColored{60}0.30 & \ValuesColored{60}0.30 & \ValuesColored{60}0.30 & \ValuesColored{58}0.29 & \ValuesColored{60}0.30 & \ValuesColored{60}0.30 \\
    IE    & \ValuesColored{44}0.22 & \ValuesColored{44}0.22 & \ValuesColored{44}0.22 & \ValuesColored{44}0.22 & \ValuesColored{44}0.22 & \ValuesColored{18}0.09 & \ValuesColored{44}0.22 & \ValuesColored{44}0.22 & \ValuesColored{44}0.22 & \ValuesColored{44}0.22 & \ValuesColored{44}0.22 & \ValuesColored{44}0.22 & \ValuesColored{44}0.22 & \ValuesColored{44}0.22 & \ValuesColored{186}0.93 & \ValuesColored{44}0.22 & \ValuesColored{44}0.22 \\
    IT    & \ValuesColored{144}0.72 & \ValuesColored{196}0.98 & \ValuesColored{200}1 & \ValuesColored{106}0.53 & \ValuesColored{166}0.83 & \ValuesColored{146}0.73 & \ValuesColored{170}0.85 & \ValuesColored{200}1 & \ValuesColored{132}0.66 & \ValuesColored{196}0.98 & \ValuesColored{116}0.58 & \ValuesColored{0}0 & \ValuesColored{196}0.98 & \ValuesColored{168}0.84 & \ValuesColored{126}0.63 & \ValuesColored{198}0.99 & \ValuesColored{124}0.62 \\
    LV    & \ValuesColored{54}0.27 & \ValuesColored{54}0.27 & \ValuesColored{54}0.27 & \ValuesColored{54}0.27 & \ValuesColored{54}0.27 & \ValuesColored{54}0.27 & \ValuesColored{54}0.27 & \ValuesColored{54}0.27 & \ValuesColored{54}0.27 & \ValuesColored{54}0.27 & \ValuesColored{54}0.27 & \ValuesColored{54}0.27 & \ValuesColored{54}0.27 & \ValuesColored{54}0.27 & \ValuesColored{54}0.27 & \ValuesColored{54}0.27 & \ValuesColored{54}0.27 \\
		LT    & \ValuesColored{80}0.40 & \ValuesColored{80}0.40 & \ValuesColored{80}0.40 & \ValuesColored{80}0.40 & \ValuesColored{80}0.40 & \ValuesColored{80}0.40 & \ValuesColored{80}0.40 & \ValuesColored{80}0.40 & \ValuesColored{80}0.40 & \ValuesColored{80}0.40 & \ValuesColored{80}0.40 & \ValuesColored{80}0.40 & \ValuesColored{80}0.40 & \ValuesColored{80}0.40 & \ValuesColored{80}0.40 & \ValuesColored{80}0.40 & \ValuesColored{80}0.40 \\
    LU    & \ValuesColored{26}0.13 & \ValuesColored{26}0.13 & \ValuesColored{26}0.13 & \ValuesColored{26}0.13 & \ValuesColored{26}0.13 & \ValuesColored{6}0.03 & \ValuesColored{26}0.13 & \ValuesColored{26}0.13 & \ValuesColored{6}0.03 & \ValuesColored{26}0.13 & \ValuesColored{48}0.24 & \ValuesColored{26}0.13 & \ValuesColored{26}0.13 & \ValuesColored{26}0.13 & \ValuesColored{66}0.33 & \ValuesColored{26}0.13 & \ValuesColored{26}0.13 \\
    NL    & \ValuesColored{46}0.23 & \ValuesColored{46}0.23 & \ValuesColored{46}0.23 & \ValuesColored{46}0.23 & \ValuesColored{46}0.23 & \ValuesColored{22}0.11 & \ValuesColored{46}0.23 & \ValuesColored{46}0.23 & \ValuesColored{36}0.18 & \ValuesColored{46}0.23 & \ValuesColored{36}0.18 & \ValuesColored{46}0.23 & \ValuesColored{46}0.23 & \ValuesColored{0}0 & \ValuesColored{80}0.40 & \ValuesColored{46}0.23 & \ValuesColored{46}0.23 \\
    PL    & \ValuesColored{40}0.20 & \ValuesColored{40}0.20 & \ValuesColored{40}0.20 & \ValuesColored{40}0.20 & \ValuesColored{40}0.20 & \ValuesColored{40}0.20 & \ValuesColored{40}0.20 & \ValuesColored{40}0.20 & \ValuesColored{40}0.20 & \ValuesColored{40}0.20 & \ValuesColored{40}0.20 & \ValuesColored{40}0.20 & \ValuesColored{40}0.20 & \ValuesColored{40}0.20 & \ValuesColored{22}0.11 & \ValuesColored{40}0.20 & \ValuesColored{40}0.20 \\
    PT    & \ValuesColored{14}0.07 & \ValuesColored{14}0.07 & \ValuesColored{14}0.07 & \ValuesColored{16}0.08 & \ValuesColored{14}0.07 & \ValuesColored{14}0.07 & \ValuesColored{14}0.07 & \ValuesColored{14}0.07 & \ValuesColored{14}0.07 & \ValuesColored{14}0.07 & \ValuesColored{14}0.07 & \ValuesColored{14}0.07 & \ValuesColored{14}0.07 & \ValuesColored{0}0 & \ValuesColored{14}0.07 & \ValuesColored{14}0.07 & \ValuesColored{14}0.07 \\
    RO    & \ValuesColored{0}0 & \ValuesColored{0}0 & \ValuesColored{0}0 & \ValuesColored{0}0 & \ValuesColored{0}0 & \ValuesColored{0}0 & \ValuesColored{0}0 & \ValuesColored{0}0 & \ValuesColored{0}0 & \ValuesColored{0}0 & \ValuesColored{0}0 & \ValuesColored{0}0 & \ValuesColored{0}0 & \ValuesColored{0}0 & \ValuesColored{0}0 & \ValuesColored{0}0 & \ValuesColored{0}0 \\
    SK    & \ValuesColored{44}0.22 & \ValuesColored{44}0.22 & \ValuesColored{44}0.22 & \ValuesColored{200}1 & \ValuesColored{44}0.22 & \ValuesColored{120}0.60 & \ValuesColored{150}0.75 & \ValuesColored{44}0.22 & \ValuesColored{74}0.37 & \ValuesColored{44}0.22 & \ValuesColored{200}1 & \ValuesColored{44}0.22 & \ValuesColored{44}0.22 & \ValuesColored{44}0.22 & \ValuesColored{26}0.13 & \ValuesColored{44}0.22 & \ValuesColored{44}0.22 \\
		SI    & \ValuesColored{114}0.57 & \ValuesColored{20}0.10 & \ValuesColored{92}0.46 & \ValuesColored{198}0.99 & \ValuesColored{114}0.57 & \ValuesColored{196}0.98 & \ValuesColored{112}0.56 & \ValuesColored{90}0.45 & \ValuesColored{196}0.98 & \ValuesColored{44}0.22 & \ValuesColored{200}1 & \ValuesColored{114}0.57 & \ValuesColored{0}0 & \ValuesColored{96}0.48 & \ValuesColored{170}0.85 & \ValuesColored{84}0.42 & \ValuesColored{114}0.57 \\
    ES    & \ValuesColored{18}0.09 & \ValuesColored{18}0.09 & \ValuesColored{0}0 & \ValuesColored{42}0.21 & \ValuesColored{18}0.09 & \ValuesColored{18}0.09 & \ValuesColored{18}0.09 & \ValuesColored{0}0 & \ValuesColored{18}0.09 & \ValuesColored{18}0.09 & \ValuesColored{18}0.09 & \ValuesColored{18}0.09 & \ValuesColored{18}0.09 & \ValuesColored{2}0.01 & \ValuesColored{0}0 & \ValuesColored{18}0.09 & \ValuesColored{18}0.09 \\
    SE    & \ValuesColored{30}0.15 & \ValuesColored{0}0 & \ValuesColored{30}0.15 & \ValuesColored{30}0.15 & \ValuesColored{30}0.15 & \ValuesColored{30}0.15 & \ValuesColored{30}0.15 & \ValuesColored{30}0.15 & \ValuesColored{30}0.15 & \ValuesColored{30}0.15 & \ValuesColored{8}0.04 & \ValuesColored{30}0.15 & \ValuesColored{30}0.15 & \ValuesColored{30}0.15 & \ValuesColored{144}0.72 & \ValuesColored{30}0.15 & \ValuesColored{30}0.15 \\
		CH    & \ValuesColored{30}0.15 & \ValuesColored{30}0.15 & \ValuesColored{30}0.15 & \ValuesColored{200}1 & \ValuesColored{30}0.15 & \ValuesColored{30}0.15 & \ValuesColored{48}0.24 & \ValuesColored{30}0.15 & \ValuesColored{0}0 & \ValuesColored{30}0.15 & \ValuesColored{10}0.05 & \ValuesColored{30}0.15 & \ValuesColored{30}0.15 & \ValuesColored{30}0.15 & \ValuesColored{52}0.26 & \ValuesColored{30}0.15 & \ValuesColored{30}0.15 \\
    GB    & \ValuesColored{34}0.17 & \ValuesColored{34}0.17 & \ValuesColored{34}0.17 & \ValuesColored{34}0.17 & \ValuesColored{34}0.17 & \ValuesColored{66}0.33 & \ValuesColored{34}0.17 & \ValuesColored{34}0.17 & \ValuesColored{0}0 & \ValuesColored{34}0.17 & \ValuesColored{2}0.01 & \ValuesColored{34}0.17 & \ValuesColored{34}0.17 & \ValuesColored{34}0.17 & \ValuesColored{80}0.40 & \ValuesColored{34}0.17 & \ValuesColored{34}0.17 \\
		\bottomrule
\end{tabu}%
\end{scriptsize}
\end{table}

\begin{table}[htbp]
\caption[]{Market power parameter values for all traders (T) in all markets (C) in April-September. The country abbreviations are given in \mbox{Table \ref{tab:P3_CountryAbbr}}. \DZ{} and \NO{} supply the \ac{EU} via the pipeline network and via \ac{LNG} shipments. This is distinguished by ${}^N$ and ${}^L$, respectively. Figures above 0.5 are colored red with increasing value, and below 0.5 green with decreasing value, to visually underpin our findings.}
\label{tab:P3_CVvaluesAP}
\centering \small \renewcommand\arraystretch{1.5}
\addtolength{\tabcolsep}{-4pt}
\begin{scriptsize}
\begin{tabu} to 0.9\linewidth{crrrrrrrrrrrrrrrrr}
    \toprule
\backslashbox{C}{T} & AZ    & DK    & DZ$^L$    & DZ$^N$   & EG    & GB    & LY    & NG    & NL    & NO$^L$    & NO$^N$   & OM    & PE    & QA    & RU    & TT    & YE \\
    \midrule
    AT    & \ValuesColored{148}0.74 & \ValuesColored{94}0.47 & \ValuesColored{104}0.52 & \ValuesColored{146}0.73 & \ValuesColored{148}0.74 & \ValuesColored{186}0.93 & \ValuesColored{170}0.85 & \ValuesColored{114}0.57 & \ValuesColored{188}0.94 & \ValuesColored{200}1 & \ValuesColored{168}0.84 & \ValuesColored{148}0.74 & \ValuesColored{200}1 & \ValuesColored{188}0.94 & \ValuesColored{144}0.72 & \ValuesColored{100}0.50 & \ValuesColored{148}0.74 \\
    BE    & \ValuesColored{68}0.34 & \ValuesColored{200}1 & \ValuesColored{80}0.40 & \ValuesColored{200}1 & \ValuesColored{68}0.34 & \ValuesColored{94}0.47 & \ValuesColored{68}0.34 & \ValuesColored{120}0.60 & \ValuesColored{32}0.16 & \ValuesColored{68}0.34 & \ValuesColored{30}0.15 & \ValuesColored{68}0.34 & \ValuesColored{68}0.34 & \ValuesColored{112}0.56 & \ValuesColored{200}1 & \ValuesColored{68}0.34 & \ValuesColored{68}0.34 \\
    BG    & \ValuesColored{122}0.61 & \ValuesColored{200}1 & \ValuesColored{200}1 & \ValuesColored{200}1 & \ValuesColored{122}0.61 & \ValuesColored{200}1 & \ValuesColored{200}1 & \ValuesColored{200}1 & \ValuesColored{200}1 & \ValuesColored{200}1 & \ValuesColored{200}1 & \ValuesColored{122}0.61 & \ValuesColored{122}0.61 & \ValuesColored{200}1 & \ValuesColored{56}0.28 & \ValuesColored{200}1 & \ValuesColored{122}0.61 \\
    HR    & \ValuesColored{188}0.94 & \ValuesColored{132}0.66 & \ValuesColored{162}0.81 & \ValuesColored{200}1 & \ValuesColored{188}0.94 & \ValuesColored{198}0.99 & \ValuesColored{196}0.98 & \ValuesColored{138}0.69 & \ValuesColored{200}1 & \ValuesColored{188}0.94 & \ValuesColored{198}0.99 & \ValuesColored{188}0.94 & \ValuesColored{188}0.94 & \ValuesColored{188}0.94 & \ValuesColored{192}0.96 & \ValuesColored{50}0.25 & \ValuesColored{188}0.94 \\
    CZ    & \ValuesColored{68}0.34 & \ValuesColored{68}0.34 & \ValuesColored{68}0.34 & \ValuesColored{28}0.14 & \ValuesColored{68}0.34 & \ValuesColored{22}0.11 & \ValuesColored{68}0.34 & \ValuesColored{68}0.34 & \ValuesColored{30}0.15 & \ValuesColored{68}0.34 & \ValuesColored{66}0.33 & \ValuesColored{68}0.34 & \ValuesColored{68}0.34 & \ValuesColored{68}0.34 & \ValuesColored{158}0.79 & \ValuesColored{68}0.34 & \ValuesColored{68}0.34 \\
    DK    & \ValuesColored{110}0.55 & \ValuesColored{72}0.36 & \ValuesColored{84}0.42 & \ValuesColored{172}0.86 & \ValuesColored{110}0.55 & \ValuesColored{88}0.44 & \ValuesColored{110}0.55 & \ValuesColored{110}0.55 & \ValuesColored{110}0.55 & \ValuesColored{110}0.55 & \ValuesColored{148}0.74 & \ValuesColored{110}0.55 & \ValuesColored{110}0.55 & \ValuesColored{28}0.14 & \ValuesColored{186}0.93 & \ValuesColored{110}0.55 & \ValuesColored{110}0.55 \\
    EE    & \ValuesColored{74}0.37 & \ValuesColored{74}0.37 & \ValuesColored{74}0.37 & \ValuesColored{74}0.37 & \ValuesColored{74}0.37 & \ValuesColored{74}0.37 & \ValuesColored{74}0.37 & \ValuesColored{74}0.37 & \ValuesColored{74}0.37 & \ValuesColored{74}0.37 & \ValuesColored{74}0.37 & \ValuesColored{74}0.37 & \ValuesColored{74}0.37 & \ValuesColored{74}0.37 & \ValuesColored{74}0.37 & \ValuesColored{74}0.37 & \ValuesColored{74}0.37 \\
    FI    & \ValuesColored{56}0.28 & \ValuesColored{56}0.28 & \ValuesColored{56}0.28 & \ValuesColored{56}0.28 & \ValuesColored{56}0.28 & \ValuesColored{56}0.28 & \ValuesColored{56}0.28 & \ValuesColored{56}0.28 & \ValuesColored{56}0.28 & \ValuesColored{56}0.28 & \ValuesColored{56}0.28 & \ValuesColored{56}0.28 & \ValuesColored{56}0.28 & \ValuesColored{56}0.28 & \ValuesColored{56}0.28 & \ValuesColored{56}0.28 & \ValuesColored{56}0.28 \\
		FR    & \ValuesColored{102}0.51 & \ValuesColored{180}0.90 & \ValuesColored{40}0.20 & \ValuesColored{194}0.97 & \ValuesColored{102}0.51 & \ValuesColored{196}0.98 & \ValuesColored{146}0.73 & \ValuesColored{188}0.94 & \ValuesColored{196}0.98 & \ValuesColored{134}0.67 & \ValuesColored{38}0.19 & \ValuesColored{102}0.51 & \ValuesColored{14}0.07 & \ValuesColored{194}0.97 & \ValuesColored{200}1 & \ValuesColored{174}0.87 & \ValuesColored{102}0.51 \\
    DE    & \ValuesColored{64}0.32 & \ValuesColored{64}0.32 & \ValuesColored{64}0.32 & \ValuesColored{114}0.57 & \ValuesColored{64}0.32 & \ValuesColored{200}1 & \ValuesColored{64}0.32 & \ValuesColored{64}0.32 & \ValuesColored{36}0.18 & \ValuesColored{64}0.32 & \ValuesColored{40}0.20 & \ValuesColored{64}0.32 & \ValuesColored{64}0.32 & \ValuesColored{200}1 & \ValuesColored{52}0.26 & \ValuesColored{64}0.32 & \ValuesColored{64}0.32 \\
    GR    & \ValuesColored{156}0.78 & \ValuesColored{156}0.78 & \ValuesColored{194}0.97 & \ValuesColored{144}0.72 & \ValuesColored{56}0.28 & \ValuesColored{126}0.63 & \ValuesColored{88}0.44 & \ValuesColored{186}0.93 & \ValuesColored{144}0.72 & \ValuesColored{126}0.63 & \ValuesColored{158}0.79 & \ValuesColored{0}0 & \ValuesColored{150}0.75 & \ValuesColored{200}1 & \ValuesColored{176}0.88 & \ValuesColored{146}0.73 & \ValuesColored{18}0.09 \\
    HU    & \ValuesColored{76}0.38 & \ValuesColored{76}0.38 & \ValuesColored{76}0.38 & \ValuesColored{200}1 & \ValuesColored{76}0.38 & \ValuesColored{164}0.82 & \ValuesColored{200}1 & \ValuesColored{76}0.38 & \ValuesColored{200}1 & \ValuesColored{76}0.38 & \ValuesColored{200}1 & \ValuesColored{76}0.38 & \ValuesColored{76}0.38 & \ValuesColored{76}0.38 & \ValuesColored{38}0.19 & \ValuesColored{76}0.38 & \ValuesColored{76}0.38 \\
    IE    & \ValuesColored{56}0.28 & \ValuesColored{56}0.28 & \ValuesColored{56}0.28 & \ValuesColored{56}0.28 & \ValuesColored{56}0.28 & \ValuesColored{20}0.10 & \ValuesColored{56}0.28 & \ValuesColored{56}0.28 & \ValuesColored{188}0.94 & \ValuesColored{56}0.28 & \ValuesColored{138}0.69 & \ValuesColored{56}0.28 & \ValuesColored{56}0.28 & \ValuesColored{56}0.28 & \ValuesColored{200}1 & \ValuesColored{56}0.28 & \ValuesColored{56}0.28 \\
    IT    & \ValuesColored{200}1 & \ValuesColored{200}1 & \ValuesColored{200}1 & \ValuesColored{200}1 & \ValuesColored{200}1 & \ValuesColored{200}1 & \ValuesColored{200}1 & \ValuesColored{200}1 & \ValuesColored{200}1 & \ValuesColored{200}1 & \ValuesColored{200}1 & \ValuesColored{200}1 & \ValuesColored{200}1 & \ValuesColored{200}1 & \ValuesColored{200}1 & \ValuesColored{200}1 & \ValuesColored{200}1 \\
		LV    & \ValuesColored{62}0.31 & \ValuesColored{62}0.31 & \ValuesColored{62}0.31 & \ValuesColored{62}0.31 & \ValuesColored{62}0.31 & \ValuesColored{62}0.31 & \ValuesColored{62}0.31 & \ValuesColored{62}0.31 & \ValuesColored{62}0.31 & \ValuesColored{62}0.31 & \ValuesColored{62}0.31 & \ValuesColored{62}0.31 & \ValuesColored{62}0.31 & \ValuesColored{62}0.31 & \ValuesColored{62}0.31 & \ValuesColored{62}0.31 & \ValuesColored{62}0.31 \\
    LT    & \ValuesColored{90}0.45 & \ValuesColored{90}0.45 & \ValuesColored{90}0.45 & \ValuesColored{90}0.45 & \ValuesColored{90}0.45 & \ValuesColored{90}0.45 & \ValuesColored{90}0.45 & \ValuesColored{90}0.45 & \ValuesColored{90}0.45 & \ValuesColored{90}0.45 & \ValuesColored{90}0.45 & \ValuesColored{90}0.45 & \ValuesColored{90}0.45 & \ValuesColored{90}0.45 & \ValuesColored{90}0.45 & \ValuesColored{90}0.45 & \ValuesColored{90}0.45 \\
    LU    & \ValuesColored{96}0.48 & \ValuesColored{44}0.22 & \ValuesColored{84}0.42 & \ValuesColored{20}0.10 & \ValuesColored{96}0.48 & \ValuesColored{146}0.73 & \ValuesColored{74}0.37 & \ValuesColored{72}0.36 & \ValuesColored{154}0.77 & \ValuesColored{96}0.48 & \ValuesColored{174}0.87 & \ValuesColored{96}0.48 & \ValuesColored{96}0.48 & \ValuesColored{66}0.33 & \ValuesColored{160}0.80 & \ValuesColored{96}0.48 & \ValuesColored{96}0.48 \\
    NL    & \ValuesColored{84}0.42 & \ValuesColored{0}0 & \ValuesColored{52}0.26 & \ValuesColored{32}0.16 & \ValuesColored{84}0.42 & \ValuesColored{66}0.33 & \ValuesColored{84}0.42 & \ValuesColored{80}0.40 & \ValuesColored{88}0.44 & \ValuesColored{84}0.42 & \ValuesColored{86}0.43 & \ValuesColored{84}0.42 & \ValuesColored{84}0.42 & \ValuesColored{46}0.23 & \ValuesColored{126}0.63 & \ValuesColored{84}0.42 & \ValuesColored{84}0.42 \\
    PL    & \ValuesColored{54}0.27 & \ValuesColored{54}0.27 & \ValuesColored{54}0.27 & \ValuesColored{54}0.27 & \ValuesColored{54}0.27 & \ValuesColored{54}0.27 & \ValuesColored{54}0.27 & \ValuesColored{54}0.27 & \ValuesColored{68}0.34 & \ValuesColored{54}0.27 & \ValuesColored{200}1 & \ValuesColored{54}0.27 & \ValuesColored{54}0.27 & \ValuesColored{54}0.27 & \ValuesColored{28}0.14 & \ValuesColored{54}0.27 & \ValuesColored{54}0.27 \\
    PT    & \ValuesColored{28}0.14 & \ValuesColored{28}0.14 & \ValuesColored{14}0.07 & \ValuesColored{60}0.30 & \ValuesColored{28}0.14 & \ValuesColored{28}0.14 & \ValuesColored{28}0.14 & \ValuesColored{0}0 & \ValuesColored{28}0.14 & \ValuesColored{28}0.14 & \ValuesColored{28}0.14 & \ValuesColored{28}0.14 & \ValuesColored{28}0.14 & \ValuesColored{12}0.06 & \ValuesColored{28}0.14 & \ValuesColored{28}0.14 & \ValuesColored{28}0.14 \\
    RO    & \ValuesColored{4}0.02 & \ValuesColored{4}0.02 & \ValuesColored{4}0.02 & \ValuesColored{4}0.02 & \ValuesColored{4}0.02 & \ValuesColored{4}0.02 & \ValuesColored{4}0.02 & \ValuesColored{4}0.02 & \ValuesColored{4}0.02 & \ValuesColored{4}0.02 & \ValuesColored{4}0.02 & \ValuesColored{4}0.02 & \ValuesColored{4}0.02 & \ValuesColored{4}0.02 & \ValuesColored{4}0.02 & \ValuesColored{4}0.02 & \ValuesColored{4}0.02 \\
    SK    & \ValuesColored{96}0.48 & \ValuesColored{200}1 & \ValuesColored{200}1 & \ValuesColored{200}1 & \ValuesColored{96}0.48 & \ValuesColored{160}0.80 & \ValuesColored{200}1 & \ValuesColored{200}1 & \ValuesColored{200}1 & \ValuesColored{96}0.48 & \ValuesColored{200}1 & \ValuesColored{96}0.48 & \ValuesColored{96}0.48 & \ValuesColored{200}1 & \ValuesColored{44}0.22 & \ValuesColored{96}0.48 & \ValuesColored{96}0.48 \\
		SI    & \ValuesColored{128}0.64 & \ValuesColored{116}0.58 & \ValuesColored{136}0.68 & \ValuesColored{186}0.93 & \ValuesColored{128}0.64 & \ValuesColored{168}0.84 & \ValuesColored{178}0.89 & \ValuesColored{120}0.60 & \ValuesColored{172}0.86 & \ValuesColored{6}0.03 & \ValuesColored{176}0.88 & \ValuesColored{128}0.64 & \ValuesColored{0}0 & \ValuesColored{146}0.73 & \ValuesColored{200}1 & \ValuesColored{82}0.41 & \ValuesColored{128}0.64 \\
    ES    & \ValuesColored{34}0.17 & \ValuesColored{34}0.17 & \ValuesColored{14}0.07 & \ValuesColored{100}0.50 & \ValuesColored{34}0.17 & \ValuesColored{34}0.17 & \ValuesColored{34}0.17 & \ValuesColored{0}0 & \ValuesColored{34}0.17 & \ValuesColored{34}0.17 & \ValuesColored{34}0.17 & \ValuesColored{34}0.17 & \ValuesColored{34}0.17 & \ValuesColored{18}0.09 & \ValuesColored{22}0.11 & \ValuesColored{34}0.17 & \ValuesColored{34}0.17 \\
    SE    & \ValuesColored{50}0.25 & \ValuesColored{6}0.03 & \ValuesColored{50}0.25 & \ValuesColored{24}0.12 & \ValuesColored{50}0.25 & \ValuesColored{32}0.16 & \ValuesColored{50}0.25 & \ValuesColored{50}0.25 & \ValuesColored{32}0.16 & \ValuesColored{50}0.25 & \ValuesColored{200}1 & \ValuesColored{50}0.25 & \ValuesColored{50}0.25 & \ValuesColored{98}0.49 & \ValuesColored{198}0.99 & \ValuesColored{50}0.25 & \ValuesColored{50}0.25 \\
    CH    & \ValuesColored{74}0.37 & \ValuesColored{74}0.37 & \ValuesColored{176}0.88 & \ValuesColored{200}1 & \ValuesColored{74}0.37 & \ValuesColored{124}0.62 & \ValuesColored{56}0.28 & \ValuesColored{74}0.37 & \ValuesColored{22}0.11 & \ValuesColored{74}0.37 & \ValuesColored{80}0.40 & \ValuesColored{74}0.37 & \ValuesColored{74}0.37 & \ValuesColored{188}0.94 & \ValuesColored{92}0.46 & \ValuesColored{74}0.37 & \ValuesColored{74}0.37 \\
		GB    & \ValuesColored{78}0.39 & \ValuesColored{78}0.39 & \ValuesColored{172}0.86 & \ValuesColored{0}0 & \ValuesColored{78}0.39 & \ValuesColored{154}0.77 & \ValuesColored{78}0.39 & \ValuesColored{78}0.39 & \ValuesColored{94}0.47 & \ValuesColored{78}0.39 & \ValuesColored{34}0.17 & \ValuesColored{78}0.39 & \ValuesColored{78}0.39 & \ValuesColored{14}0.07 & \ValuesColored{158}0.79 & \ValuesColored{78}0.39 & \ValuesColored{78}0.39 \\
    \bottomrule
\end{tabu}%
\end{scriptsize}
\end{table}

\begin{table}[htbp]
\caption{Country abbreviations. The attribution of a country into the ``Eastern \ac{EU}'' or ``Western \ac{EU}'' is determined by its geographical location and does not coincide with the historical political division of Europe.}
\label{tab:P3_CountryAbbr}
\centering \small \renewcommand\arraystretch{1.5}
\begin{tabu} to 0.9\linewidth{lXlXlX}
		\multicolumn{2}{l}{\textbf{Eastern \ac{EU}}} 	& \multicolumn{2}{l}{\textbf{Western \ac{EU}}} 	& \multicolumn{2}{l}{\textbf{Non-EU}} \\
		\multicolumn{2}{l}{\textbf{consumers}} 				& \multicolumn{2}{l}{\textbf{consumers}} 				& \multicolumn{2}{l}{\textbf{suppliers}} \\
			BG& \BG{}& AT& \AT{}& AZ& \AZ{}\\ 
			CZ& \CZ{}& BE& \BE{}& DZ$^L$& \DZ{} (LNG)\\ 
			EE& \EE{}& CH& \CH{}& DZ$^N$& \DZ{} (Pipeline)\\
			FI& \FI{}& DE& \DE{}& EG& \EG{}\\
			GR& \GR{}& DK& \DK{}& LY& \LY{}\\ 
			HR& \HR{}& ES& \ES{}& NG& Nigeria\\
			HU& \HU{}& FR& \FR{}& NO$^L$& \NO{} (LNG)\\
			LT& \LT{}& GB& \UK{}& NO$^N$& \NO{} (Pipeline)\\
			LV& \LV{}& IE& \IE{}& OM& \OM{}\\
			PL& \PL{}& IT& \IT{}& PE& \PE{}\\ 
			RO& \RO{}& LU& \LU{}& QA& \QA{}\\
			SI& \SI{}& NL& \TNL{}& RU& \RU\\
			SK& \SK{}& PT& \PT{}& TT& \TT{}\\ 
			& 		 & SE& \SE{}& YE& \YE{}\\  
\end{tabu}
\end{table}
\mbox{Tables \ref{tab:P3_CVvaluesOC}} \mbox{and \ref{tab:P3_CVvaluesAP}} show that the market power parameters $\MP$ are similar per country and not per trader, for instance, \IT{} faces high $\MP$, whereas \RO{} faces very low $\MP$. This originates from the fact that the inverse demand function is equal for all traders, whereas the marginal supply costs $\phi$ vary per trader \textit{and} country. The highest mean values were found for the large pipeline-bound suppliers \RU{} (0.50), \DZ{}$^N$ (0.46), \NO{}$^N$ (0.44), \UK{} (0.43), \LY{} (0.41), and \tNL{} (0.40), which mirrors the observations we make in reality. 

Moreover, the $\MP$'s are generally \textit{higher} in the European summer than in winter, which is rather counter-intuitive from an economic point of view; one would expect traders to exert less market power in times demand and prices are low. However, this finding can be explained by the interplay of the market power parameters $\MP$ and $\lambdaO$. On one hand, Equation \eqref{eqn:P3_MPCalibChoice} implies that for non-zero reference sales $\qRef_f$ the market power $\MP_f$ of a trader $f$ increases with increasing difference between the market price and its marginal supply costs $\lambdaO-\phi_f$. On the other hand, if $\lambdaORange \cup \lambdaRefRange$ contains more than one point, and $\qRef$, $\sRef$, $\phiRef$, and $\etaRef$ are given and fixed, the choice of $\lambdaO$ determines the $\MP_f$, and any combination of $\lambdaO \in \lambdaRefRange$ and $\MP_f(\lambdaO,\qRef_f,\sRef,\phiRef_f,\etaRef) \in [0,1]$ is admissible. This renders a certain flexibility when choosing the parameters and the corresponding outcome. 
In our example, we chose reference prices $\lambdaRef$ to be equal per country in the summer and winter periods. This corresponds to fixing a rather low $\lambdaO \in \lambdaORange \cup \lambdaRefRange$ for the winter period, and a rather high $\lambdaO \in \lambdaORange \cup \lambdaRefRange$ for the summer period. As a consequence, our calibration gives rather low $\MP$'s in winter and rather high $\MP$'s in summer.
We carried out additional simulations, for which the results are not shown here, with a spread of \unit[10]{\%} and \unit[20]{\%} between summer and winter reference prices. For the former setting, the average $\MP$ is similar in both seasons, while in the latter the average $\MP$ is clearly higher in the high demand period. We conclude that the reference data greatly influences the obtained equilibrium and the corresponding parameters, and therefore should be carefully chosen.

\clearpage
\section{Model equations and exemplary market setting}  \label{app:P3_ModelEq_Pic_Notation}
\mbox{\ref{app:P3_ModelEq_Pic_Notation}} is largely reproduced from \citet{Baltensperger2015}. We introduce the model equations in \mbox{\ref{app:P3_ModelEquationsSubSec}}, show an exemplary market setting with two interconnected nodes in \mbox{\ref{app:P3_figure}}, and introduce the associated notation in \mbox{\ref{app:P3_NotationSubSec}}. Note that we follow the convention used by \citet{Baltensperger2015} and include producers in the notion of service providers, although producers are not an infrastructure service.

\subsection{Model equations} \label{app:P3_ModelEquationsSubSec}
Equations \eqref{eqn:P3_NoLOSSnoTauNoS} describe the mechanics of the spatial partial equilibrium model of the European gas market in detail. We refrain from showing the loss terms in the equations to achieve a more compact notation. 
\begin{subequations} \label{eqn:P3_NoLOSSnoTauNoS}
\allowdisplaybreaks
\begin{alignat}{4}
0 &\leq& \, \LINC^P_{nt} + \QUAC^P_{nt} q^P_{fnt} + \alpha^P_{nt} + \alpha^{PT}_{n} - \phi^N_{fnt} & \perp q^P_{fnt} &&\,\geq 0 \quad \forall f,n,t \label{eqn:P3_dqP}\\
0 &\leq&\, \LINC^I_{nt} + \alpha^I_{nt} + \alpha^{IT}_{n} + \phi^N_{fnt} - \phi^S_{fn}& \perp q^I_{fnt} &&\,\geq 0 \quad \forall f,n,t \label{eqn:P3_dqI}\\
0 &\leq&\, \LINC^X_{zt} + \alpha^X_{nt} + \alpha^{XT}_{n} - \phi^N_{fnt} + \phi^S_{fn} & \perp q^X_{fnt} &&\,\geq 0 \quad \forall f,n,t \label{eqn:P3_dqX}\\
0 &\leq&\, \LINC^A_{nmt} + \alpha^A_{nmt} + \alpha^{AT}_{nm} - \phi^N_{fmt} + \phi^N_{fnt} & \perp q^A_{fnmt} &&\,\geq 0 \quad \forall f,n,m,t \label{eqn:P3_dqA}\\
0 &\leq&\, \LINC^L_{nt} + \alpha^L_{nt} + \alpha^{LT}_{n} + \LINC^B_{nmt} + \alpha^B_{nmt} + \alpha^{BT}_{nm} &&& \notag\\
&&\, + \LINC^R_{mt} + \alpha^R_{mt} + \alpha^{RT}_{m} - \phi^N_{fmt} +  \phi^N_{fnt} & \perp q^B_{fnmt} &&\,\geq 0 \quad \forall f,n,m,t \label{eqn:P3_dqB}\\
0 &\leq&\, -\lambda^C_{nt} -\theta^C_{fnt} \SLP^C_{nt} q^C_{fnt} + \phi^N_{fnt}& \perp q^C_{fnt} &&\,\geq 0  \quad \forall f,n,t\\
0 &\leq&\, q^{P}_{fnt} + q^{X}_{fnt} + \sum \limits_{m \in \mathcal{A}(n)} q_{fmnt}^{A} + \sum \limits_{m \in \mathcal{B}(n)} q_{fmnt}^{B}&&& \notag \\
&&\,- q^{I}_{fnt} -  q^{C}_{fnt} - \sum \limits_{m \in \mathcal{A}(n)} q_{fnmt}^{A} - \sum \limits_{m \in \mathcal{B}(n)} q_{fnmt}^{B} & \perp \phi^N_{fnt} &&\,\geq 0  \quad \forall f,n,t \label{eqn:P3_dphiN}\\
0 &\leq&\, \sum \limits_{t \in \mathcal{T}} q^{I}_{fnt} - \sum \limits_{t \in \mathcal{T}} q^{X}_{fnt} &\perp \phi^S_{fn} &&\,\geq 0  \quad \forall f,n\label{eqn:P3_dphiS}\\
0 &\leq&\, \overline{\CAP}^Z_{zt} - \sum \limits_{f \in \mathcal{F}(z)} q^{Z}_{fzt} & \perp \alpha^Z_{zt} &&\,\geq 0  \quad \forall z,t \label{eqn:P3_dalphaZ}\\
0 &\leq&\, \overline{\CAP}^{ZT}_{z} - \sum \limits_{t \in \mathcal{T}} \sum \limits_{f \in \mathcal{F}(z)} q^{Z}_{fzt} & \perp \alpha^{ZT}_z &&\,\geq 0  \quad \forall z \label{eqn:P3_dbetaZ}\\
0 &\leq&\, \lambda^C_{nt} - \left( \INT^C_{nt} + \SLP^C_{nt} \sum \limits_{f \in \mathcal{F}(n)} q^{C}_{fnt} \right) &\perp \lambda^C_{nt} &&\,\geq 0  \quad \forall n,t
\end{alignat}
\end{subequations}

\clearpage
\subsection{Graphical illustration} \label{app:P3_figure}
\begin{figure}[!htbp]
\centering
\includegraphics[width=0.9\textwidth]{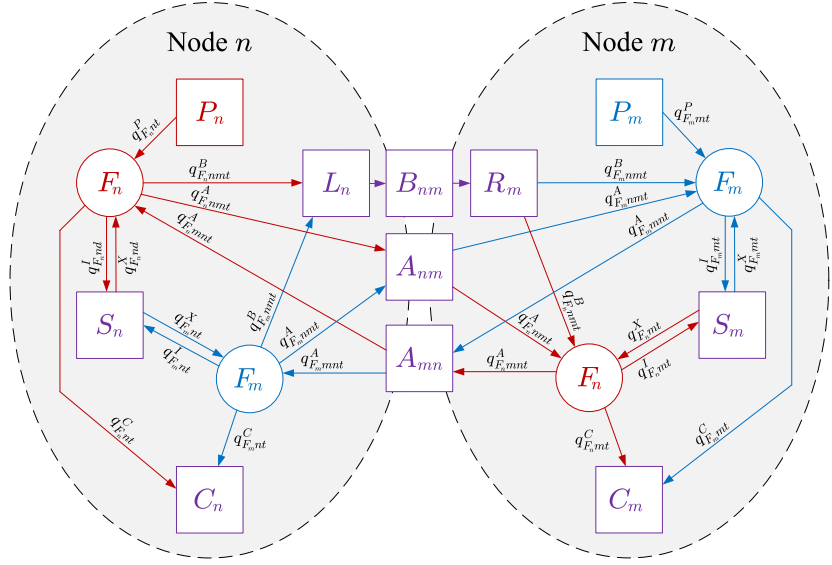}
\caption{Gas market model with two nodes. $P$: producer. $A$: pipeline operator. $L$: liquefaction plant operator. $B$: LNG shipment. $R$: regasification plant operator. $S$: storage operator. $C$: consumer. $F_n$: trader associated with producer $P_n$. $F_m$: trader associated with producer $P_m$. $q^{P}_{fnt}$: quantity delivered from producer to trader $f$ at node $n$ in time period $t$. $q^{A}_{fnmt}$: pipeline transportation of trader $f$ via arc $nm$ in period $t$. $q^{B}_{fnmt}$: LNG shipment of trader $f$ via arc $nm$ in period $t$. $q^{I}_{fnt}$: storage injection by trader $f$ at node $n$ in period $t$. $q^{X}_{fnt}$: storage extraction by trader $f$ at node $n$ in period $t$. $q^{C}_{fnt}$: sales by trader $f$ to consumer in node $n$ in period $t$. The traders $F_n$ and $F_m$, there decision variables, and their corresponding producers $P_n$ and $P_m$ are colored red ($n$) and blue ($m$), respectively. Service providers, except producers, and flows between them are marked purple, as well as the consumers, since all traders trade with them.}
\label{fig:P3_Model_full}
\end{figure}

\clearpage
\subsection{Notation} \label{app:P3_NotationSubSec}
\begin{table}[htbp]
\caption{This table introduces the nomenclature concerning service providers, traders and consumers.}
\label{tab:P3_ServiceProviders}
\centering \small \renewcommand\arraystretch{1.5}
\begin{tabu} to 0.9\linewidth{lX}
\toprule 
\multicolumn{2}{l}{\textbf{Service providers, traders and consumers}} \\
\midrule 
$A_{nm}$& Transmission system operator of \mbox{pipeline $nm$}\\
$B_{nm}$& Shipping company transporting \ac{LNG} from $n$ \mbox{to $m$}\\
$C_n$& Consumer at \mbox{node $n$}\\
$I_n$& Storage operator injecting gas at \mbox{node $n$}\\
$L_n$& Liquefaction plant operator at \mbox{node $n$}\\
$P_n$& Gas producing company at \mbox{node $n$}\\
$R_n$& Regasification plant operator at \mbox{node $n$}\\
$S_n$& Storage operator at \mbox{node $n$}\\
$F_n$& The trader associated with producer at \mbox{node $n$}\\
$X_n$& Storage operator extracting gas at \mbox{node $n$}\\
$Z_z$& Placeholder for a service provider ($P_n$, $I_n$, $X_n$, $L_n$, $R_n$, $A_{nm}$, $B_{nm}$) at node $n$ / \mbox{arc $nm$}\\
\bottomrule
\end{tabu}
\end{table}

{\small \singlespacing \renewcommand\arraystretch{1.5}
\begin{longtabu}  to 0.9\linewidth{lX}
\caption{This table introduces all sets used for the mathematical description of the model.}
\label{tab:P3_Sets} \\
\toprule 
\multicolumn{2}{l}{\textbf{Sets}}  \\
\midrule 
\endfirsthead
\multicolumn{2}{c}{\begin{footnotesize}\tablename\ \thetable\ -- \textit{Continued from previous page}\end{footnotesize}} \\
\toprule
\multicolumn{2}{l}{\textbf{Sets}}  \\
\midrule
\endhead
\bottomrule
\multicolumn{2}{r}{\begin{footnotesize}\textit{Continued on next page}\end{footnotesize}} \\
\endfoot
\bottomrule
\endlastfoot
$t \in \mathcal{T} = \{T_1, \ldots, T_{\bar{t}}\}$& A time period $t$ in the set $\mathcal{T}$ of all periods of a year\\
$n, m \in \mathcal{N} = \{N_1, \ldots, N_{\bar{n}}\}$& Nodes $n,m$ in the set $\mathcal{N}$ of all nodes \\ 
$f \in \mathcal{F} = \{F_1, \ldots, F_{\bar{n}}\}$& A trader $f$ in the set $\mathcal{F}$ of all traders\\
$z \in \mathcal{Z}$ & A node/arc element from the \mbox{set $\mathcal{Z}$} \\
\midrule
$\mathcal{A} \subset \mathcal{N} \times \mathcal{N}$& Set of arcs connecting $2$ nodes by pipeline					\\
$\mathcal{B} \subset \mathcal{N} \times \mathcal{N}$ & Set of arcs connecting $2$ nodes by ship					\\	
$\mathcal{C} \subseteq \mathcal{N}$ & Set of nodes at which a consumer is active \\
$\mathcal{I} \subseteq \mathcal{N}$ & Set of nodes at which storage injection	is possible					\\
$\mathcal{L} \subseteq \mathcal{N}$ & Set of nodes at which a liquefaction terminal operator is active				\\	
$\mathcal{P} \subseteq \mathcal{N}$ & Set of nodes at which a gas producer is active	\\							
$\mathcal{R} \subseteq \mathcal{N}$ & Set of nodes at which a regasification terminal operator is active		\\
$\mathcal{X} \subseteq \mathcal{N}$ & Set of nodes at which storage extraction is possible						\\
$\mathcal{Z} \in \{\mathcal{P},\mathcal{L},\mathcal{B},\mathcal{R},\mathcal{A},\mathcal{I},\mathcal{X}\}$ & Placeholder for the set of nodes/arcs at which a type of service provider is active \\
\midrule
$\mathcal{A}(n) \subseteq \mathcal{N} \setminus \{n\}$ & Set of nodes which are connected to $n$ by pipeline \\
$\mathcal{B}(n) \subseteq \mathcal{N} \setminus \{n\}$ & Set of nodes which are connected to $n$ by ship \\
$\mathcal{C}(f) \subseteq \mathcal{N}$& The set of all nodes with consumers which are reachable by \mbox{trader $f$}\\
$\mathcal{N}(f) \subseteq \mathcal{N}$& The set of all nodes which are reachable by \mbox{trader $f$}\\
$\mathcal{F}(z)$& The set of all traders active at node/\mbox{arc $z$}\\
$\mathcal{Z}(f)$& The set of all nodes/arcs in which service $Z$ is active and are reachable by \mbox{trader $f$}\\ 
\end{longtabu}
}

{\small \singlespacing \renewcommand\arraystretch{1.5}
\begin{longtabu}  to 0.9\linewidth{lX}
\caption{The \textit{parameters} are generally described by capital Roman letters. Lower-case Roman letters are only chosen if the parameter is directly linked to a variable of the same name. Occasionally, lower-case Greek letters are chosen to follow conventions. The superscripts indicate whether the parameter is related to a service provider of type $Z \in \{P,L,B,R,A,I,X\}$ or a consumer $C$. Subscripts indicate the trader $f$, node/arc $z$, and the period of the year $t$ the parameter is related to.}
\label{tab:P3_Parameters} \\
\toprule 
\multicolumn{2}{l}{\textbf{Parameters}}  \\
\midrule 
\endfirsthead
\multicolumn{2}{c}{\begin{footnotesize}\tablename\ \thetable\ -- \textit{Continued from previous page}\end{footnotesize}} \\
\toprule
\multicolumn{2}{l}{\textbf{Parameters}}  \\
\midrule
\endhead
\bottomrule
\multicolumn{2}{r}{{\footnotesize\textit{Continued on next page}}} \\
\endfoot
\bottomrule
\endlastfoot
$\overline{\CAP}^Z_{nt}$		& Maximum capacity of service $Z$ located at $z$ in \mbox{period $t$}\\
$\overline{\CAP}^{ZT}_{n}$	& Maximum capacity of service $Z$ located at $z$ over all \mbox{periods $\mathcal{T}$}\\
$\INT^C_{nt}$ 							& Maximum willingness to pay (intercept of inverse demand function with the $s^C_{nt}=0$ - axis) of consumers at node $n$ in \mbox{period $t$} \\
$\LINC^Z_{zt}$ 							& Linear cost function term for service $Z$ located at $z$ in \mbox{period	$t$}\\
$\LOSS^Z_z $								& Loss factor when using service $Z$  located \mbox{at $z$} \\
$\QUAC^Z_{zt}$ 							& Quadratic cost function term for service $Z$ located at $z$ in \mbox{period $t$}\\
$\SLP^C_{nt}$								& Slope of the inverse demand curve of the consumers at node $n$ in \mbox{period $t$}, is assumed strictly negative\\
$\MP_{fnt}$									& Market power parameter of trader $f$ in node $n$ and \mbox{period $t$}\\ 
\end{longtabu}
}

\begin{table}[htbp]
\caption{The \textit{variables} are described by lowercase letters. Primal variables are Roman, while dual variables are Greek letters. The superscripts indicate whether the variable is related to a service provider of type $Z \in \{P,L,B,R,A,I,X\}$, a consumer $C$, or a node $N$. Subscripts indicate the trader $f$ the variable corresponds to, at which node/arc $z$ the transaction or service is located, and in which period of the year $t$ it takes place.}
\label{tab:P3_Variables}
\centering \small \renewcommand\arraystretch{1.5}
\begin{tabu} to 0.9\linewidth{lX}
\toprule
\multicolumn{2}{l}{\textbf{Variables}} \\
\midrule
$q^{C}_{fnt}$ & Flow of trader $f$ to consumer $C$ at node $n$ in \mbox{period $t$}\\
$q^{Z}_{fzt}$ & Flow between trader $f$ and service provider $Z$ at node/arc $z$ in \mbox{period $t$} \\
$\alpha^Z_{zt}$ & Congestion fee of service $Z$ located at $z$ in \mbox{period $t$}\\
$\alpha^{ZT}_{z}$ & Congestion fee on annual usage of service $Z$ located at $z$ \\
$\phi^N_{fnt}$ & Dual variable of the volume balance of trader $f$ at node $n$ and \mbox{period $t$}\\
$\phi^S_{fn}$ & Dual variable of the annual volume balance of trader $f$ in storage $S$ at \mbox{node $n$}\\
$\lambda^C_{nt}$ & Wholesale price at node $n$ in \mbox{period $t$}\\
\bottomrule
\end{tabu}
\end{table}

\begin{table}[htbp]
\caption{This table introduces the \textit{functions}. The superscript $C$ indicates that the function is related to the wholesale market. Subscripts indicate at which node $n$ and in which period of the year $t$ the function is valid.}
\label{tab:P3_Functions}
\centering \small \renewcommand\arraystretch{1.5}
\begin{tabu} to 0.9\linewidth{lX}
\toprule 
\multicolumn{2}{l}{\textbf{Functions}} \\
\midrule 
$\Lambda^{C}_{nt}(s^C_{nt})$ & Inverse demand function of consumer $C$ at node $n$ in \mbox{period $t$}.\\
\bottomrule
\end{tabu}
\end{table}

\clearpage

\clearpage
\bibliography{00_Wrapper}

\begin{thebibliography}{23}
\expandafter\ifx\csname natexlab\endcsname\relax\def\natexlab#1{#1}\fi
\providecommand{\url}[1]{\texttt{#1}}
\providecommand{\href}[2]{#2}
\providecommand{\path}[1]{#1}
\providecommand{\DOIprefix}{doi:}
\providecommand{\ArXivprefix}{arXiv:}
\providecommand{\URLprefix}{URL: }
\providecommand{\Pubmedprefix}{pmid:}
\providecommand{\doi}[1]{\href{http://dx.doi.org/#1}{\path{#1}}}
\providecommand{\Pubmed}[1]{\href{pmid:#1}{\path{#1}}}
\providecommand{\bibinfo}[2]{#2}
\ifx\xfnm\relax \def\xfnm[#1]{\unskip,\space#1}\fi
\bibitem[{Baltensperger et~al.(2015)Baltensperger, F{\"{u}}chslin, Kr{\"{u}}tli
  \& Lygeros}]{Baltensperger2015}
\bibinfo{author}{Baltensperger, T.}, \bibinfo{author}{F{\"{u}}chslin, R.~M.},
  \bibinfo{author}{Kr{\"{u}}tli, P.}, \& \bibinfo{author}{Lygeros, J.}
  \bibinfo{year}{2015}.
\newblock \bibinfo{title}{{Multiplicity of equilibria in conjectural variation
  models of natural gas markets}}.
\newblock \href{http://arxiv.org/abs/1510.04473}{\tt arXiv:1510.04473}.
\bibitem[{Chyong \& Hobbs(2014)}]{Chyong2014e}
\bibinfo{author}{Chyong, C.~K.}, \& \bibinfo{author}{Hobbs, B.~F.}
  \bibinfo{year}{2014}.
\newblock \bibinfo{title}{{Strategic Eurasian natural gas market model for
  energy security and policy analysis: Formulation and application to South
  Stream: Online supplementary material}}.
\newblock {\it \bibinfo{journal}{Energy Economics}\/},  {\it
  \bibinfo{volume}{44}\/}, \bibinfo{pages}{198--211}.
\bibitem[{Cobanli(2014)}]{Cobanli2014}
\bibinfo{author}{Cobanli, O.} \bibinfo{year}{2014}.
\newblock \bibinfo{title}{{Central Asian gas in Eurasian power game}}.
\newblock {\it \bibinfo{journal}{Energy Policy}\/},  {\it
  \bibinfo{volume}{68}\/}, \bibinfo{pages}{348--370}.
  \DOIprefix\doi{10.1016/j.enpol.2013.12.027}.
\bibitem[{D{\'{\i}}az et~al.(2014)D{\'{\i}}az, Campos \& Villar}]{Diaz2014}
\bibinfo{author}{D{\'{\i}}az, C.~a.}, \bibinfo{author}{Campos, F.~A.}, \&
  \bibinfo{author}{Villar, J.} \bibinfo{year}{2014}.
\newblock \bibinfo{title}{{Existence and uniqueness of conjectured supply
  function equilibria}}.
\newblock {\it \bibinfo{journal}{International Journal of Electrical Power and
  Energy Systems}\/},  {\it \bibinfo{volume}{58}\/}, \bibinfo{pages}{266--273}.
  \DOIprefix\doi{10.1016/j.ijepes.2014.01.027}.
\bibitem[{D{\'{\i}}az et~al.(2010)D{\'{\i}}az, Villar, Campos \&
  Reneses}]{Diaz2010}
\bibinfo{author}{D{\'{\i}}az, C.~a.}, \bibinfo{author}{Villar, J.},
  \bibinfo{author}{Campos, F.~A.}, \& \bibinfo{author}{Reneses, J.}
  \bibinfo{year}{2010}.
\newblock \bibinfo{title}{{Electricity market equilibrium based on conjectural
  variations}}.
\newblock {\it \bibinfo{journal}{Electric Power Systems Research}\/},  {\it
  \bibinfo{volume}{80}\/}, \bibinfo{pages}{1572--1579}.
  \DOIprefix\doi{10.1016/j.epsr.2010.07.012}.
\bibitem[{D{\'{\i}}az et~al.(2011)D{\'{\i}}az, Villar, Campos \&
  Rodr{\'{\i}}guez}]{Diaz2011}
\bibinfo{author}{D{\'{\i}}az, C.~a.}, \bibinfo{author}{Villar, J.},
  \bibinfo{author}{Campos, F.~A.}, \& \bibinfo{author}{Rodr{\'{\i}}guez,
  M.~{\'{A}}.} \bibinfo{year}{2011}.
\newblock \bibinfo{title}{{A new algorithm to compute conjectured supply
  function equilibrium in electricity markets}}.
\newblock {\it \bibinfo{journal}{Electric Power Systems Research}\/},  {\it
  \bibinfo{volume}{81}\/}, \bibinfo{pages}{384--392}.
  \DOIprefix\doi{10.1016/j.epsr.2010.10.002}.
\bibitem[{{\noopsort{EC}}{\ac{EC}}(2013{\natexlab{a}})}]{EuropeanCommission2013a}
\bibinfo{author}{{\noopsort{EC}}{\ac{EC}}} \bibinfo{year}{2013}{\natexlab{a}}.
\newblock {\it \bibinfo{title}{{Quarterly report on European Gas Markets, first
  quarter 2013}}\/}.
\newblock \bibinfo{type}{Technical Report} \bibinfo{address}{retrieved from
  \url{https://ec.europa.eu/energy/en/statistics/market-analysis}}.
\bibitem[{{\noopsort{EC}}{\ac{EC}}(2013{\natexlab{b}})}]{EuropeanCommission2013}
\bibinfo{author}{{\noopsort{EC}}{\ac{EC}}} \bibinfo{year}{2013}{\natexlab{b}}.
\newblock {\it \bibinfo{title}{{Quarterly report on European Gas Markets,
  second quarter 2013}}\/}.
\newblock \bibinfo{type}{Technical Report} \bibinfo{address}{retrieved from
  \url{https://ec.europa.eu/energy/en/statistics/market-analysis}}.
\bibitem[{{\noopsort{EC}}{\ac{EC}}(2015)}]{EuropeanCommission}
\bibinfo{author}{{\noopsort{EC}}{\ac{EC}}} \bibinfo{year}{2015}.
\newblock \bibinfo{title}{{Eurostat}}.
\newblock
  \bibinfo{howpublished}{\url{http://ec.europa.eu/eurostat/de/data/database}}.
\newblock \bibinfo{annote}{{[Online; accessed 13-July-2015]}}.
\bibitem[{{\noopsort{ENTSO-G}}{\ac{ENTSO-G}}(2012)}]{ENTSO-G2012c}
\bibinfo{author}{{\noopsort{ENTSO-G}}{\ac{ENTSO-G}}} \bibinfo{year}{2012}.
\newblock {\it \bibinfo{title}{{Ten-Year Network Development Plan}}\/}.
\newblock \bibinfo{type}{Technical Report} \bibinfo{address}{retrieved from
  \url{http://www.entsog.eu/publications/tyndp}}.
\bibitem[{Garc{\'{\i}}a-Alcalde et~al.(2002)Garc{\'{\i}}a-Alcalde, Ventosa,
  Rivier, Ramos \& Rela{\~{n}}o}]{Garcia-Alcalde2002}
\bibinfo{author}{Garc{\'{\i}}a-Alcalde, A.}, \bibinfo{author}{Ventosa, M.},
  \bibinfo{author}{Rivier, M.}, \bibinfo{author}{Ramos, A.}, \&
  \bibinfo{author}{Rela{\~{n}}o, G.} \bibinfo{year}{2002}.
\newblock \bibinfo{title}{{Fitting electricity market models. A conjectural
  variations approach.}}
\newblock In {\it \bibinfo{booktitle}{14th PSCC, Sevilla}\/}.
\bibitem[{Huppmann \& Egging(2014)}]{Huppmann2014}
\bibinfo{author}{Huppmann, D.}, \& \bibinfo{author}{Egging, R.}
  \bibinfo{year}{2014}.
\newblock \bibinfo{title}{{Market power, fuel substitution and infrastructure
  – A large-scale equilibrium model of global energy markets}}.
\newblock {\it \bibinfo{journal}{Energy}\/},  {\it \bibinfo{volume}{75}\/},
  \bibinfo{pages}{483--500}. \DOIprefix\doi{10.1016/j.energy.2014.08.004}.
\bibitem[{{\noopsort{IEA}}{\ac{IEA}}(2012)}]{InternationalEnergyAgencyIEA2012c}
\bibinfo{author}{{\noopsort{IEA}}{\ac{IEA}}} \bibinfo{year}{2012}.
\newblock {\it \bibinfo{title}{{IEA statistics: Natural gas information}}\/}.
\newblock \bibinfo{type}{Technical Report}.
\bibitem[{{\noopsort{IEA}}{\ac{IEA}}(2013{\natexlab{a}})}]{InternationalEnergyAgencyIEA2013}
\bibinfo{author}{{\noopsort{IEA}}{\ac{IEA}}}
  \bibinfo{year}{2013}{\natexlab{a}}.
\newblock {\it \bibinfo{title}{{IEA statistics: Natural gas information}}\/}.
\newblock \bibinfo{type}{Technical Report}.
\bibitem[{{\noopsort{IEA}}{\ac{IEA}}(2013{\natexlab{b}})}]{Iea2013}
\bibinfo{author}{{\noopsort{IEA}}{\ac{IEA}}}
  \bibinfo{year}{2013}{\natexlab{b}}.
\newblock {\it \bibinfo{title}{{World energy outlook 2013}}\/}.
\newblock \bibinfo{type}{Technical Report}.
\newblock \DOIprefix\doi{10.1787/weo-2013-en}.
\bibitem[{{\noopsort{IGU}}{\ac{IGU}}(2014)}]{Union2014}
\bibinfo{author}{{\noopsort{IGU}}{\ac{IGU}}} \bibinfo{year}{2014}.
\newblock {\it \bibinfo{title}{{Wholesale Gas Price Survey \textendash 2014
  Edition}}\/}.
\newblock \bibinfo{type}{Technical Report} \bibinfo{address}{retrieved from
  \url{http://www.igu.org/research/wholesale-gas-price-survey-2014-edition}}.
\bibitem[{Kamiński(2011)}]{Kaminski2011}
\bibinfo{author}{Kamiński, J.} \bibinfo{year}{2011}.
\newblock \bibinfo{title}{{Market power in a coal-based power generation
  sector: The case of Poland}}.
\newblock {\it \bibinfo{journal}{Energy}\/},  {\it \bibinfo{volume}{36}\/},
  \bibinfo{pages}{6634--6644}. \DOIprefix\doi{10.1016/j.energy.2011.08.048}.
\bibitem[{Lise et~al.(2008)Lise, Hobbs \& van Oostvoorn}]{Lise2008}
\bibinfo{author}{Lise, W.}, \bibinfo{author}{Hobbs, B.~F.}, \&
  \bibinfo{author}{van Oostvoorn, F.} \bibinfo{year}{2008}.
\newblock \bibinfo{title}{{Natural gas corridors between the EU and its main
  suppliers: Simulation results with the dynamic GASTALE model}}.
\newblock {\it \bibinfo{journal}{Energy Policy}\/},  {\it
  \bibinfo{volume}{36}\/}, \bibinfo{pages}{1890--1906}.
  \DOIprefix\doi{10.1016/j.enpol.2008.01.042}.
\bibitem[{Liu et~al.(2006)Liu, Lie \& Lo}]{Liu2006}
\bibinfo{author}{Liu, J.}, \bibinfo{author}{Lie, T.}, \& \bibinfo{author}{Lo,
  K.} \bibinfo{year}{2006}.
\newblock \bibinfo{title}{{An Empirical Method of Dynamic Oligopoly Behavior
  Analysis in Electricity Markets}}.
\newblock {\it \bibinfo{journal}{IEEE Transactions on Power Systems}\/},  {\it
  \bibinfo{volume}{21}\/}, \bibinfo{pages}{499--506}.
  \DOIprefix\doi{10.1109/TPWRS.2006.873054}.
\bibitem[{Liu et~al.(2007)Liu, Ni, Wu \& Cai}]{Liu2007}
\bibinfo{author}{Liu, Y.}, \bibinfo{author}{Ni, Y.~X.}, \bibinfo{author}{Wu,
  F.~F.}, \& \bibinfo{author}{Cai, B.} \bibinfo{year}{2007}.
\newblock \bibinfo{title}{{Existence and uniqueness of consistent conjectural
  variation equilibrium in electricity markets}}.
\newblock {\it \bibinfo{journal}{International Journal of Electrical Power and
  Energy Systems}\/},  {\it \bibinfo{volume}{29}\/}, \bibinfo{pages}{455--461}.
  \DOIprefix\doi{10.1016/j.ijepes.2006.11.006}.
\bibitem[{{L{\'{o}}pez de Haro} et~al.(2007){L{\'{o}}pez de Haro},
  {S{\'{a}}nchez Mart{\'{\i}}n}, {de la Hoz Ardiz} \& {Fern{\'{a}}ndez
  Caro}}]{LopezdeHaro2007}
\bibinfo{author}{{L{\'{o}}pez de Haro}, S.}, \bibinfo{author}{{S{\'{a}}nchez
  Mart{\'{\i}}n}, P.}, \bibinfo{author}{{de la Hoz Ardiz}, J.~E.}, \&
  \bibinfo{author}{{Fern{\'{a}}ndez Caro}, J.} \bibinfo{year}{2007}.
\newblock \bibinfo{title}{{Estimating conjectural variations for electricity
  market models}}.
\newblock {\it \bibinfo{journal}{European Journal of Operational Research}\/},
  {\it \bibinfo{volume}{181}\/}, \bibinfo{pages}{1322--1338}.
  \DOIprefix\doi{10.1016/j.ejor.2005.12.039}.
\bibitem[{Tremblay \& Tremblay(2012)}]{Tremblay2012}
\bibinfo{author}{Tremblay, V.~J.}, \& \bibinfo{author}{Tremblay, C.~H.}
  \bibinfo{year}{2012}.
\newblock {\it \bibinfo{title}{{New Perspectives on Industrial Organization:
  With Contributions from Behavioral Economics and Game Theory}}\/}.
\newblock \bibinfo{address}{New York}: \bibinfo{publisher}{Springer}.
\newblock \DOIprefix\doi{10.1007/978-1-4614-3241-8}.
\bibitem[{{\noopsort{UN}}{\ac{UN}}(2009)}]{UnitedNations}
\bibinfo{author}{{\noopsort{UN}}{\ac{UN}}} \bibinfo{year}{2009}.
\newblock \bibinfo{title}{{UN data -- a world of information}}.
\newblock \bibinfo{howpublished}{\url{http://data.un.org/}}.
\newblock \bibinfo{annote}{{[Online; accessed 09-October-2014]}}.

\end{thebibliography}
\end{document}